\documentclass[12pt]{amsart}

\usepackage{amsmath}
\usepackage{amscd}
\usepackage{amssymb}
\usepackage{latexsym}
\usepackage{stmaryrd}

\usepackage[arrow,curve,matrix,tips,2cell]{xy}

\newtheorem{theorem}{Theorem}[section]
\newtheorem*{theorem*}{Theorem}
\newtheorem{lemma}[theorem]{Lemma}
\newtheorem*{lemma*}{Lemma}
\newtheorem{corollary}[theorem]{Corollary}
\newtheorem*{corollary*}{Corollary}
\newtheorem{proposition}[theorem]{Proposition}

\newtheorem{remark}[theorem]{Remark}
\newtheorem{question}[theorem]{Question}
\newtheorem{definition}[theorem]{Definition}

\newtheorem{example}[theorem]{Example}
\newtheorem{examples}[theorem]{Examples}
\newcommand{\bgl}{\begin{equation}} 
\newcommand{\egl}{\end{equation}}
\newcommand{\bgloz}{\begin{equation*}} 
\newcommand{\egloz}{\end{equation*}}
\newcommand{\bgln}{\begin{eqnarray}} 
\newcommand{\egln}{\end{eqnarray}}
\newcommand{\bglnoz}{\begin{eqnarray*}} 
\newcommand{\eglnoz}{\end{eqnarray*}}
\newcommand{\btheo}{\begin{theorem}}
\newcommand{\etheo}{\end{theorem}}
\newcommand{\btheooz}{\begin{theorem*}}
\newcommand{\etheooz}{\end{theorem*}}
\newcommand{\blemma}{\begin{lemma}}
\newcommand{\elemma}{\end{lemma}}
\newcommand{\blemmaoz}{\begin{lemma*}}
\newcommand{\elemmaoz}{\end{lemma*}}
\newcommand{\bproof}{\begin{proof}}
\newcommand{\eproof}{\end{proof}}
\newcommand{\bbew}{\begin{beweis}}
\newcommand{\ebew}{\end{beweis}}
\newcommand{\bremark}{\begin{remark}\em}
\newcommand{\eremark}{\end{remark}}
\newcommand{\bquestion}{\begin{question}\em}
\newcommand{\equestion}{\end{question}}
\newcommand{\bdefin}{\begin{definition}}
\newcommand{\edefin}{\end{definition}}
\newcommand{\bprop}{\begin{proposition}}
\newcommand{\eprop}{\end{proposition}}
\newcommand{\bcor}{\begin{corollary}}
\newcommand{\ecor}{\end{corollary}}
\newcommand{\bcoroz}{\begin{corollary*}}
\newcommand{\ecoroz}{\end{corollary*}}
\newcommand{\bfa}{\begin{cases}} 
\newcommand{\efa}{\end{cases}}
\newcommand{\bexample}{\begin{example}\em}
\newcommand{\eexample}{\end{example}}
\newcommand{\bexamples}{\begin{examples}\em}
\newcommand{\eexamples}{\end{examples}}
\newcommand{\cB}{\mathcal B}
\newcommand{\cC}{\mathcal C}

\newcommand{\cE}{\mathcal E}
\newcommand{\cF}{\mathcal F}
\newcommand{\cG}{\mathcal G}

\newcommand{\cJ}{\mathcal J}

\newcommand{\cL}{\mathcal L}

\newcommand{\cO}{\mathcal O}
\newcommand{\cP}{\mathcal P}

\newcommand{\cR}{\mathcal R}
\newcommand{\cS}{\mathcal S}

\newcommand{\cZ}{\mathcal Z}
\def\Cz{\mathbb{C}}
\def\Fz{\mathbb{F}}

\def\Nz{\mathbb{N}}

\def\Zz{\mathbb{Z}}

\def\1z{\mathbb{1}}

\newcommand{\an}[1]{``#1''} 
\newcommand{\ti}{\tilde}

\newcommand{\lori}{\longrightarrow}

\newcommand{\ma}{\mapsto} 
\newcommand\onto{\twoheadrightarrow} 
\newcommand\into{\hookrightarrow} 
\newcommand{\Rarr}{\Rightarrow} 
\newcommand{\Larr}{\Leftarrow} 
\newcommand{\LRarr}{\Leftrightarrow} 

\newcommand{\ve}{\varepsilon}

\def\SEMI{\mbox{$\times\kern-2pt\vrule height5pt width.6pt \kern3pt $}}

\newcommand{\img}{{\rm Im\,}}

\newcommand{\Spec}{{\rm Spec\,}} 

\newcommand{\id}{{\rm id}}

\newcommand{\reg}{^\times} 
\newcommand{\lspan}{{\rm span}} 
\newcommand{\clspan}{\overline{\lspan}} 
\newcommand{\abs}[1]{\left|#1\right|} 
\newcommand{\norm}[1]{\left\|#1\right\|} 
\newcommand{\defeq}{\mathrel{:=}} 

\newcommand{\dop}{\text{: }} 

\newcommand{\dom}{{\rm dom}}

\newcommand{\lge}{\left\{} 
\newcommand{\rge}{\right\}} 
\newcommand{\lru}{\left(} 
\newcommand{\rru}{\right)} 
\newcommand{\lsp}{\left\langle} 
\newcommand{\rsp}{\right\rangle} 
\newcommand{\rukl}[1]{\lru #1 \rru} 
\newcommand{\gekl}[1]{\lge #1 \rge} 
\newcommand{\spkl}[1]{\lsp #1 \rsp} 

\newcommand{\tfin}{\text{fin}} 
\newcommand{\tinf}{\text{inf}}

\newcommand{\menge}[2]{\gekl{ #1 \dop #2 }} 
\title{Partial transformation groupoids attached to graphs and semigroups}

\author{Xin Li}

\address{Xin Li, School of Mathematical Sciences, Queen Mary University of London, Mile End Road, London E1 4NS, United Kingdom}
\email{xin.li@qmul.ac.uk}


\begin{document}

\begin{abstract}
We introduce the notion of continuous orbit equivalence for partial dynamical systems, and give an equivalent characterization in terms of Cartan-isomorphisms for partial C*-crossed products. Both graph C*-algebras and semigroup C*-algebras can be described as C*-algebras attached to partial dynamical systems. As applications, for graphs, we generalize and explain a result of Matsumoto and Matui relating orbit equivalence and Cartan-isomorphism, and for semigroups, we strengthen several structural results for semigroup C*-algebras concerning amenability, nuclearity as well as simplicity of boundary quotients. We also discuss pure infiniteness for partial transformation groupoids arising from graphs and semigroups.
\end{abstract}

\thanks{Research supported by EPSRC grant EP/M009718/1.}

\maketitle


\setlength{\parindent}{0cm} \setlength{\parskip}{0.5cm}

\section{Introduction}

Recently, in the setting of ordinary topological dynamical systems, the notion of continuous orbit equivalence was introduced, and a C*-algebraic characterization was given involving Cartan-isomorphisms (see \cite{Li4}). We present a generalization to partial dynamical systems. This step is important as many C*-algebras appear naturally as crossed products attached to partial dynamical systems, while the setting of ordinary dynamical systems seems rather restricted. Our main motivation stems from graph C*-algebras and semigroup C*-algebras, both of which can be described in a very natural way as C*-algebras of partial dynamical systems. This description turns out to be very helpful for the study of structural properties of these C*-algebras.

In the case of graphs, we obtain a very easy explanation why all graph C*-algebras are nuclear. This goes back to the observation that non-abelian free semigroups embed into amenable groups. Furthermore, we generalize and explain results in \cite{MM} about orbit equivalence and Cartan isomorphism from shifts of finite type to general graphs. We note that a generalization of the results in \cite{MM} has been established independently in \cite{BCW}.

For semigroup C*-algebras and their boundary quotients, we are able to generalize several structural results from \cite{Li2}. The key idea is that using partial transformation groupoids, we obtain structural results for our C*-algebras without assuming independence or the Toeplitz condition, which were crucial in our previous approach (see \cite{Li1,Li2}). We obtain general characterizations for nuclearity of semigroup C*-algebras and simplicity of boundary quotients.

Furthermore, we study which partial transformation groupoids of graphs and semigroups are purely infinite, in the sense of \cite{Ma2}. In the case of graphs, we are able to prove that the partial system of a graph is residually topologically free and purely infinite if and only if the corresponding graph C*-algebra is purely infinite. For semigroups, we show that the partial transformation groupoid corresponding to a boundary quotient is purely infinite as long as the semigroup is not left reversible. We are also able to identify a class of integral domains whose $ax+b$-semigroups have purely infinite partial transformation groupoids. These results strengthen and explain previous results in \cite{Li2,Li3}.

%

\section{Partial actions, transformation groupoids, C*-algebras and Cartan subalgebras}

In the following, groups are discrete and countable, and topological spaces are locally compact, Hausdorff and second countable. 

\setlength{\parindent}{0cm} \setlength{\parskip}{0cm}
\bdefin
Let $G$ be a group with identity $e$, and let $X$ be a topological space. A partial action $\alpha$ of $G$ on $X$ consists of
\begin{itemize}
\item a collection $\gekl{U_g}_{g \in G}$ of open subsets $U_g \subseteq X$,
\item a collection $\gekl{\alpha_g}_{g \in G}$ of homeomorphisms $\alpha_g: \: U_{g^{-1}} \to U_g, \, x \ma g.x$ such that
\begin{itemize}
\item $U_e = X$, $\alpha_e = \id_X$;
\item for all $g_1, g_2 \in G$, we have $g_2.(U_{(g_1 g_2)^{-1}} \cap U_{g_2^{-1}}) = U_{g_2} \cap U_{g_1^{-1}}$, and $(g_1 g_2).x = g_1.(g_2.x)$ for all $x \in U_{(g_1 g_2)^{-1}} \cap U_{g_2^{-1}}$.
\end{itemize}
\end{itemize}
We call such a triple $(X,G,\alpha)$ a partial system, and denote it by $\alpha: \: G \curvearrowright X$ or simply $G \curvearrowright X$. 
\edefin
\setlength{\parindent}{0cm} \setlength{\parskip}{0cm}

Let $\alpha: \: G \curvearrowright X$ be a partial system. The dual action $\alpha^*$ of $\alpha$ is the partial action (in the sense of \cite{McCl}) of $G$ on $C_0(X)$ given by $\alpha^*_g: \: C_0(U_{g^{-1}}) \to C_0(U_g), \, f \ma f(g^{-1}.\sqcup)$.
\setlength{\parindent}{0cm} \setlength{\parskip}{0.25cm}

The transformation groupoid attached to the partial system $\alpha: \: G \curvearrowright X$ is given by
$$G \mathbin{_{\alpha} \ltimes} X \defeq \menge{(g,x) \in G \times X}{g \in G, x \in U_{g^{-1}}},$$
with source map $s(g,x) = x$, range map $r(g,x) = g.x$, composition $(g_1,g_2.x)(g_2,x) = (g_1g_2,x)$ and inverse $(g,x)^{-1} = (g^{-1},g.x)$. We equip $G \mathbin{_{\alpha} \ltimes} X$ with the subspace topology from $G \times X$. Usually, we write $G \ltimes X$ for $G \mathbin{_{\alpha} \ltimes} X$ if the action $\alpha$ is understood. The unit space of $G \ltimes X$ coincides with $X$. Since $G$ is discrete, $G \ltimes X$ is an \'{e}tale groupoid. Actually, if we set $G_x \defeq \menge{g \in G}{x \in U_{g^{-1}}}$ and $G^x \defeq \menge{g \in G}{x \in U_g}$ for $x \in X$, then we have canonical identifications $s^{-1}(x) \cong G_x, \, (g,x) \ma g$ and $r^{-1}(x) \cong G^x, \, (g,g^{-1}.x) \ma g$.
\setlength{\parindent}{0cm} \setlength{\parskip}{0.5cm}

Let us now recall the construction (from \cite{McCl}) of the reduced crossed product $C_0(X) \rtimes_{\alpha^*,r} G$ attached to our partial system $\alpha: \: G \curvearrowright X$. As with groupoids, we omit $\alpha^*$ in our notation for the crossed product. First of all $C_0(X) \rtimes^{\ell^1} G \defeq \menge{\sum_g f_g \delta_g \in \ell^1(G,C_0(X))}{f_g \in C_0(U_g)}$ becomes a *-algebra under component-wise addition, multiplication given by $\rukl{\sum_g f_g \delta_g} \cdot \rukl{\sum_h \ti{f}_h \delta_h} \defeq \sum_{g,h} \alpha^*_g(\alpha^*_{g^{-1}}(f_g) \ti{f}_h) \delta_{gh}$ and involution $\rukl{\sum_g f_g \delta_g}^* \defeq \sum_g \alpha^*_g(f_{g^{-1}}^*) \delta_g$.
\setlength{\parindent}{0.5cm} \setlength{\parskip}{0cm}

As in \cite{McCl}, we construct a representation of $C_0(X) \rtimes^{\ell^1} G$. Viewing $X$ as a discrete set, we define $\ell^2 X$ and the representation $M: \: C_0(X) \to \cL(\ell^2 X), \, f \ma M(f)$, where $M(f)$ is the multiplication operator $M(f)(\xi) \defeq f \cdot \xi$ for $\xi \in \ell^2 X$. $M$ is obviously a faithful representation of $C_0(X)$. Every $g \in G$ leads to a twist of $M$, namely $M_g: \: C_0(X) \to \cL(\ell^2 X)$ given by $M_g(f) \xi \defeq f \vert_{U_g}(g. \sqcup) \cdot \xi \vert_{U_{g^{-1}}}$. Here we view $f \vert_{U_g}(g. \sqcup)$ as an element in $C_b(U_{g^{-1}})$, and $C_b(U_{g^{-1}})$ acts on $\ell^2 U_{g^{-1}}$ just by multiplication operators. Given $\xi \in \ell^2 X$, we set $\xi \vert_{U_{g^{-1}}}(x) \defeq \xi(x)$ if $x \in U_{g^{-1}}$ and $\xi \vert_{U_{g^{-1}}}(x) \defeq 0$ if $x \notin U_{g^{-1}}$. In other words, $\xi \vert_{U_{g^{-1}}}$ is the component of $\xi$ in $\ell^2 U_{g^{-1}}$ with respect to the decomposition $\ell^2 X = \ell^2 U_{g^{-1}} \oplus \ell^2 U_{g^{-1}}^c$. So we have $M_g(f) \xi (x) = f(g.x) \xi(x)$ if $x \in U_{g^{-1}}$ and $M_g(f) \xi (x) = 0$ if $x \notin U_{g^{-1}}$.

Consider now the Hilbert space $H \defeq \ell^2(G, \ell^2 X) \cong \ell^2 G \otimes \ell^2 X$, and define the representation $\mu: \: C_0(X) \to \cL(H)$ given by $\mu(f)(\delta_g \otimes \xi) \defeq \delta_g \otimes M_g(f) \xi$. For $g \in G$, let $E_g$ be the orthogonal projection onto $\overline{\mu(C_0(U_{g^{-1}})) H}$. Moreover, let $\lambda$ denote the left regular representation of $G$ on $\ell^2 G$, and set $V_g \defeq (\lambda_g \otimes I) \cdot E_g$. Here $I$ is the identity operator on $H$.

We can now define the representation $\mu \times \lambda: \: C_0(X) \rtimes^{\ell^1} G \to \cL(H), \, \sum_g f_g \delta_g \ma \sum_g \mu(f_g) V_g$. Following the original definition in \cite{McCl}, we set $C_0(X) \rtimes_r G \defeq \overline{C_0(X) \rtimes^{\ell^1} G}^{\norm{\cdot}_{\mu \times \lambda}}$.
\setlength{\parindent}{0cm} \setlength{\parskip}{0.5cm}

The following result follows from \cite[Theorem~3.3]{Aba}, but we give a direct and short proof.
\setlength{\parindent}{0cm} \setlength{\parskip}{0cm}
\bprop
The canonical homomorphism
\bgl
\label{cG->}
  C_c(G \ltimes X) \to C_0(X) \rtimes^{\ell^1} G, \, \theta \ma \sum_g \theta(g,g^{-1}.\sqcup) \delta_g,
\egl
where $\theta(g,g^{-1}.\sqcup)$ is the function $U_{g^{-1}} \to \Cz, \, x \ma \theta(g,g^{-1}.x)$, extends to an isomorphism $C^*_r(G \ltimes X) \overset{\cong}{\lori} C_0(X) \rtimes_r G$.
\eprop
\bproof
\setlength{\parindent}{0.5cm} \setlength{\parskip}{0cm}
As above, let $\mu \times \lambda$ be the representation $C_0(X) \rtimes^{\ell^1} G \to \cL(H)$ which we used to define $C_0(X) \rtimes_r G$. Our first observation is
\bgl
\label{nondeg}
  \overline{\img(\mu \times \lambda)(H)}
  = \bigoplus_{h \in G} \delta_h \otimes \ell^2 U_{h^{-1}}.
\egl
To see this, observe that for all $g \in G$, $\img(E_g) \subseteq \bigoplus_h \delta_h \otimes \ell^2(U_{h^{-1}} \cap U_{(gh)^{-1}})$. This holds since for $x \notin h^{-1}.(U_h \cap U_{g^{-1}}) = U_{(gh)^{-1}} \cap U_{h^{-1}}$, $f \vert_{U_h} (h.x) = 0$ for $f \in C_0(U_{g^{-1}})$. Therefore, $\pi(C_0(U_{g^{-1}}))(\delta_h \otimes \ell^2 X) \subseteq \delta_h \otimes \ell^2(U_{h^{-1}} \cap U_{(gh)^{-1}})$. Hence $\img(E_g) \subseteq \bigoplus_h \delta_h \otimes \ell^2(U_{h^{-1}} \cap U_{(gh)^{-1}})$, and thus, $\img(V_g) \subseteq \bigoplus_h \delta_{gh} \otimes \ell^2(U_{h^{-1}} \cap U_{(gh)^{-1}}) \subseteq \bigoplus_h \delta_h \otimes \ell^2 U_{h^{-1}}$. This shows \an{$\subseteq$} in \eqref{nondeg}. For \an{$\supseteq$}, note that for $f \in C_0(X)$, $(\mu \times \lambda)(f \delta_e) = \mu(f) E_e$, and for $\xi \in \ell^2 U_{h^{-1}}$, $\mu(f) E_e (\delta_h \otimes \xi) = \delta_h \otimes f \vert_{U_h}(h. \sqcup) \xi$. So $(\mu \times \lambda)(f \delta_e)(H)$ contains $\delta_h \otimes f \cdot \xi$ for all $f \in C_0(U_{h^{-1}})$ and $\xi \in \ell^2 U_{h^{-1}}$, hence also $\delta_h \otimes \ell^2 U_{h^{-1}}$. This proves \an{$\supseteq$}.

For $x \in X$, let $G_x = \menge{g \in G}{x \in U_{g^{-1}}}$ as before. Our second observation is that for every $x \in X$, the subspace $H_x \defeq \ell^2 G_x \otimes \delta_x$ is $(\mu \times \lambda)$-invariant. It is clear that $\mu(f)$ leaves $H_x$ invariant for all $f \in C_0(X)$. For $g, h \in G$, $E_g(\delta_h \otimes \delta_x) = \delta_h \otimes \delta_x$ if $x \in U_{h^{-1}} \cap U_{(gh)^{-1}}$, and if that is the case, then $V_g(\delta_h \otimes \delta_x) = \delta_{gh} \otimes \delta_x \in H_x$.

Therefore, $H = \rukl{\bigoplus_{x \in X} H_x} \oplus (\mu \times \lambda)(C_0(X) \rtimes^{\ell^1} G)(H)^{\perp}$ is a decomposition of $H$ into $\mu \times \lambda$-invariant subspaces. For $x \in X$, set $\rho_x \defeq (\mu \times \lambda) \vert_{H_x}$. Then $C_0(X) \rtimes_r G = \overline{C_0(X) \rtimes^{\ell^1} G}^{\norm{\cdot}_{\bigoplus_x \rho_x}}$.

Moreover, we have for $x \in U_{h^{-1}}$,
\bgln
  \rho_x \rukl{\sum_g f_g \delta_g}(\delta_h \otimes \delta_x)
  &=& \sum_g \mu(f_g) V_g (\delta_h \otimes \delta_x)
  = \sum_{g: \: x \in U_{(gh)^{-1}}} \mu(f_g) (\delta_{gh} \otimes \delta_x) \nonumber \\
\label{rho_x}
  &=& \sum_{g: \: x \in U_{(gh)^{-1}}} \delta_{gh} \otimes f_g(gh.x) \delta_x
  = \sum_{k \in G_x} \delta_k \otimes f_{kh^{-1}}(k.x) \delta_x 
\egln

Let us now compare this construction with the construction of the reduced groupoid C*-algebra of $G \ltimes X$. Obviously, \eqref{cG->} is an embedding of $C_c(G \ltimes X)$ as a subalgebra which is $\norm{\cdot}_{\ell^1}$-dense in $C_0(X) \rtimes^{\ell^1} G$. Therefore, $C_0(X) \rtimes_r G = \overline{C_c(G \ltimes X)}^{\norm{\cdot}_{\bigoplus_x \rho_x}}$.

Now, to construct the reduced groupoid C*-algebra $C^*_r(G \ltimes X)$, we follow \cite[\S~2.3.4]{R09} and construct for every $x \in X$ the representation $\pi_x: \: C_c(G \ltimes X) \to \cL(\ell^2(s^{-1}(x)))$ by setting $\pi_x(\theta)(\xi)(\zeta) \defeq \sum_{\eta \, \in \, s^{-1}(x)} \theta(\zeta \eta^{-1}) \xi(\eta)$. In our case, using $s^{-1}(x) = G_x \times \gekl{x}$, we obtain for $\xi = \delta_h \otimes \delta_x$ with $h \in G_x$: $\pi_x(\theta)(\delta_h \otimes \delta_x)(k,x) = \theta((k.x)(h,x)^{-1}) = \theta(kh^{-1},h.x)$. Thus,
\bgl
\label{pi_x}
  \pi_x(\theta)(\delta_h \otimes \delta_x)(k,x) = \sum_{k \in G_x} \theta(kh^{-1},h.x) \delta_k \otimes \delta_x.
\egl
By definition, $C^*_r(G \ltimes X) = \overline{C_c(G \ltimes X)}^{\norm{\cdot}_{\bigoplus_x \pi_x}}$. Therefore, in order to show that $\norm{\cdot}_{\bigoplus_x \rho_x}$ and $\norm{\cdot}_{\bigoplus_x \rho_x}$ coincide on $C_c(G \ltimes X)$, it suffices to show that for every $x \in X$, $\pi_x$ and the restriction of $\rho_x$ to $C_c(G \ltimes X)$ are unitarily equivalent. Given $x \in X$, using $s^{-1}(x) = G_x \times \gekl{x}$, we obtain the canonical unitary $\ell^2(s^{-1}(x)) \cong H_x = \ell^2(G_x) \otimes \delta_x$, so that we may think of both $\rho_x$ and $\pi_x$ as representations on $\ell^2(G_x) \otimes \delta_x$. We then have for $x \in X$, $\theta \in C_c(G \ltimes X)$ and $h \in G_x$:
\bgloz
  \rho_x(\theta)(\delta_h \otimes \delta_x)
  \overset{\eqref{cG->}}{=}
  \rho_x(\sum_g \theta(g,g^{-1}. \sqcup) \delta_g)(\delta_h \otimes \delta_x)
  \overset{\eqref{rho_x}}{=}
  \sum_{k \in G_x} \delta_k \otimes \theta(kh^{-1},h.x) \delta_x \overset{\eqref{pi_x}}{=} \pi_x(\theta)(\delta_h \otimes \delta_x).
\egloz
This yields the canonical identification $C_0(X) \rtimes_r G \cong C^*_r(G \ltimes X)$, as desired.
\setlength{\parindent}{0cm} \setlength{\parskip}{0cm}
\eproof
Following \cite{ELQ}, we define topological freeness as follows:
\bdefin
A partial system $G \curvearrowright X$ is called topologically free if for every $e \neq g \in G$, $\menge{x \in U_{g^{-1}}}{g.x \neq x}$ is dense in $U_{g^{-1}}$.
\edefin
\setlength{\parindent}{0cm} \setlength{\parskip}{0.5cm}

\blemma
\label{TF_intersection}
A partial system $G \curvearrowright X$ is topologically free if and only if
\setlength{\parindent}{0cm} \setlength{\parskip}{0cm}

$\menge{x \in X}{g.x \neq x \ for \ all \ e \neq g \in G_x}$ is dense in $X$.
\elemma
\setlength{\parindent}{0.5cm} \setlength{\parskip}{0cm}
\bproof
The direction \an{$\Larr$} is simple: If $\menge{x \in X}{g.x \neq x \ for \ all \ e \neq g \in G_x}$ is dense in $X$, then in particular, for every fixed $e \neq g \in G$, the set $\menge{x \in X}{g.x \neq x \ if \ g \in G_x}$ is dense in $X$. Hence $\menge{x \in U_{g^{-1}}}{g.x \neq x} = \menge{x \in X}{g.x \neq x \ if \ g \in G_x} \cap U_{g^{-1}}$ is dense in $U_{g^{-1}}$.
\setlength{\parindent}{0.5cm} \setlength{\parskip}{0cm}

For \an{$\Rarr$}, note that because of topological freeness, we know that for every $e \neq g \in G$, the open set $\menge{x \in U_{g^{-1}}}{g.x \neq x} \cup \overline{U_{g^{-1}}}^c$ is dense in $X$. Therefore,
$$
  \menge{x \in X}{g.x \neq x \ for \ all \ e \neq g \in G_x}
  = \bigcap_{e \neq g \in G} \menge{x \in X}{g.x \neq x \ if \ g \in G_x}
$$
must be dense in $X$ since it contains $\bigcap_{e \neq g \in G} \rukl{ \menge{x \in U_{g^{-1}}}{g.x \neq x} \cup \overline{U_{g^{-1}}}^c }$ which is dense in $X$ by the Baire category theorem.
\eproof
\setlength{\parindent}{0cm} \setlength{\parskip}{0.5cm}

\bcor
\label{TF-->tp}
A partial system $G \curvearrowright X$ is topologically free if and only if the transformation groupoid $G \ltimes X$ is topologically principal.
\ecor
\setlength{\parindent}{0.5cm} \setlength{\parskip}{0cm}
\bproof
Recall (see \cite{R08}) that a topological groupoid $\cG$ is called topologically principal if the set of points in $\cG^{(0)}$ with trivial isotropy is dense in $\cG^{(0)}$. Here, $x \in \cG^{(0)}$ is said to have trivial isotropy if for all $\gamma \in \cG$, $s(\gamma) = t(\gamma) = x$ already implies $\gamma = x$. In the case of the transformation groupoid $\cG = G \ltimes X$, $x \in X$ has trivial isotropy if whenever $g \in G$ satisfies $g \in G_x$ and $g.x = x$, then we must have $g = e$. Hence the set of points with trivial isotropy is nothing else but $\menge{x \in X}{g.x \neq x \ for \ all \ e \neq g \in G_x}$. With this observation, our corollary follows immediately from Lemma~\ref{TF_intersection}.
\eproof
\setlength{\parindent}{0cm} \setlength{\parskip}{0.25cm}

Let us now introduce the notion of continuous orbit equivalence for partial systems, generalizing \cite[Definition~2.5]{Li4}.
\bdefin
\label{COE}
Partial systems $G \curvearrowright X$ and $H \curvearrowright Y$ are called continuously orbit equivalent if there exists a homeomorphism $\varphi: \: X \overset{\cong}{\lori} Y$ and continuous maps $a: \: \bigcup_{g \in G} \gekl{g} \times U_{g^{-1}} \to H$, $b: \: \bigcup_{h \in H} \gekl{h} \times V_{h^{-1}} \to G$ (where $V_{h^{-1}}$ is the domain of the partial homeomorphism attached to $h \in H$) such that
\bgln
\label{varphi}
  \varphi(g.x) &=& a(g,x).\varphi(x), \\
\label{phi}
  \varphi^{-1}(h.y) &=& b(h,y).\varphi^{-1}(y).
\egln
Implicitly, we require here that $a(g,x) \in H_{\varphi(x)}$ and $b(h,y) \in G_{\varphi^{-1}(y)}$.
\edefin
Note that in particular, $\varphi(G_x.x) = H_{\varphi(x)}.\varphi(x)$.

In analogy to \cite[Theorem~1.2]{Li4}, we obtain
\setlength{\parindent}{0.5cm} \setlength{\parskip}{0cm}
\btheo
\label{COE-gpd-C}
Let $G \curvearrowright X$ and $H \curvearrowright Y$ be topologically free partial systems. Then the following are equivalent:
\begin{enumerate}
\item[(i)] $G \curvearrowright X$ and $H \curvearrowright Y$ are continuously orbit equivalent,
\item[(ii)] the transformation groupoids $G \ltimes X$ and $H \ltimes Y$ are isomorphic as topological groupoids,
\item[(iii)] there exists an isomorphism $\Phi: \: C_0(X) \rtimes_r G \overset{\cong}{\lori} C_0(Y) \rtimes_r H$ with $\Phi(C_0(X)) = C_0(Y)$.
\end{enumerate}

Moreover, \an{(ii) $\Rarr$ (i)} holds in general (i.e., without the assumption of topological freeness).
\etheo
\bproof
The proof is completely analogous to the one of \cite[Theorem~1.2]{Li4}. Therefore, we refer the reader to \cite{Li4} for details, and only mention the key ideas.
\setlength{\parindent}{0.5cm} \setlength{\parskip}{0cm}

For \an{(i) $\Rarr$ (ii)}, observe that $G \ltimes X \to H \ltimes Y, \, (g,x) \ma (a(g,x),\varphi(x))$ and $H \ltimes Y \to G \ltimes X, \, (h,y) \ma (b(h,y),\varphi^{-1}(y))$ are continuous homomorphisms of groupoids which are inverse to each other. Here $\varphi$, $a$ and $b$ are as in Definition~\ref{COE}. This uses topological freeness as in \cite{Li4}.

For \an{(ii) $\Rarr$ (i)}, let $\chi: \: G \ltimes X \to H \ltimes Y$ be an isomorphism of topological groupoids. Set $\varphi \defeq \chi \vert_X$, $a$ as the composition $\bigcup_{g \in G} \gekl{g} \times U_{g^{-1}} \to G \ltimes X \overset{\chi}{\lori} H \ltimes Y \to H$, and $b$ as the composition $\bigcup_{h \in H} \gekl{h} \times V_{h^{-1}} \to H \ltimes Y \overset{\chi^{-1}}{\lori} G \ltimes X \to G$. Then it is easy to check that $\varphi$, $a$ and $b$ satisfy the conditions in Definition~\ref{COE}. This does not use topological freeness.

For \an{(ii) $\LRarr$ (iii)}, observe that by Corollary~\ref{TF-->tp} and \cite[Theorem~5.2]{R08}, $(C_0(X) \rtimes_r G,C_0(X)) \cong (C^*_r(G \ltimes X),C_0(X))$ and $(C_0(Y) \rtimes_r H,C_0(Y)) \cong (C^*_r(H \ltimes Y),C_0(Y))$ are Cartan pairs, in the sense of \cite{R08}. Then apply \cite[Proposition~4.13]{R08}.
\eproof
\setlength{\parindent}{0cm} \setlength{\parskip}{0.5cm}

\section{Examples: Inverse semigroups, graphs, and subsemigroups of groups}

\subsection{Inverse semigroups}
\label{S}

Let $S$ be an inverse semigroup with zero $0$ and $E$ the semilattice of idempotents of $S$. Assume that $\sigma$ is a partial homomorphism from $S$ to a group $G$ which is idempotent pure, i.e., $\sigma$ is a map $S\reg \to G$, where $S\reg = S \setminus \gekl{0}$, such that $\sigma(st) = \sigma(s) \sigma(t)$ for all $s,t \in S\reg$ with $st \neq 0$, and that $\sigma^{-1}(e) = E\reg (= E \setminus \gekl{0})$. Existence of such an idempotent pure partial homomorphism is equivalent to saying that $S$ is strongly $0$-$E$-unitary.
\setlength{\parindent}{0.5cm} \setlength{\parskip}{0cm}

In this situation, we describe a partial action $G \curvearrowright C^*(E)$ such that the (left) reduced C*-algebra $C^*_{\lambda}(S)$ of $S$ is canonically isomorphic to $C^*(E) \rtimes_r G$. This is the analogue of \cite[Theorem~5.2]{MS}, but without the assumption that $G$ is the universal group. Our observation gives us more freedom in the choice of idempotent pure partial homomorphisms. This will lead to a very simple criterion for nuclearity of inverse semigroup C*-algebras.
\setlength{\parindent}{0cm} \setlength{\parskip}{0.25cm}

First, recall the definition of $C^*_{\lambda}(S)$: For $s \in S$, let $\lambda_s: \: \ell^2 S\reg \to \ell^2 S\reg, \, \delta_x \ma \delta_{sx}$ if $s^*s \geq xx^*$ and $\delta_x \ma 0$ else. The reduced C*-algebra $C^*_{\lambda}(S)$ of $S$ is the sub-C*-algebra of $\cL(\ell^2 S\reg)$ generated by $\menge{\lambda_s}{s \in S}$. Note that by construction, $\lambda_0 = 0$. This is why we work with $\ell^2 S\reg$ instead of $\ell^2 S$. Let $C^*(E) \defeq C^*(\menge{\lambda_e}{e \in E}) \subseteq C^*_{\lambda}(S)$. In the following, we write $e$ for $\lambda_e \in C^*(E)$. The full inverse semigroup C*-algebra $C^*(S)$ of an inverse semigroup $S$ is the universal C*-algebra for *-representations of $S$ by partial isometries. We mod out $0$ if $S$ has a zero, as in \cite{Nor}.

Now let us describe the partial action $G \curvearrowright C^*(E)$. For $g \in G$, let $D_g$ be the sub-C*-algebra (actually ideal) of $C^*(E)$ given by $D_{g^{-1}} = \clspan(\menge{s^*s}{s \in S\reg, \sigma(s) = g})$. As $\sigma$ is idempotent pure, we have $D_e = C^*(E)$. For every $g \in G$, we have a C*-isomorphism $\alpha^*_g: \: D_{g^{-1}} \to D_g, \, s^*s \ma ss^*$. The corresponding dual action is given as follows: Let $\widehat{E} = \Spec(C^*(E))$, and for every $g \in G$, set $U_g = \Spec(D_g) \subseteq \widehat{E}$. It is easy to see that $U_{g^{-1}} = \menge{\chi \in \widehat{E}}{\chi(s^*s) = 1 \ {\rm for} \ {\rm some} \ s \in S\reg \ {\rm with} \ \sigma(s) = g}$. $\alpha_g: \: U_{g^{-1}} \to U_g$ is given by $\alpha_g(\chi) = \chi \circ \alpha^*_{g^{-1}}$. Given $\chi \in U_{g^{-1}}$ and $s \in S\reg$ with $\sigma(s) = g$ and $\chi(s^*s) = 1$, we have $\alpha_g(\chi)(e) = \chi(s^*es)$.

The same proof as in \cite{MS} gives an explicit identification of the universal groupoid of $S$ with the partial transformation groupoid $G \ltimes \widehat{E}$, in the sense of \cite[\S~4.3]{Pat}. Hence we obtain
\bprop
\label{S-EG}
We have $C^*_{\lambda}(S) \cong C^*_r(G \ltimes \widehat{E})$ and $C^*(S) \cong C^*(G \ltimes \widehat{E})$, and there are isomorphisms $C^*_{\lambda}(S) \to C^*(E) \rtimes_r G, \, \lambda_s \ma (ss^*) u_{\sigma(s)}$ and $C^*(S) \to C^*(E) \rtimes G, \, s \ma (ss^*) \delta_{\sigma(s)}$.
\eprop
 
In many situations, one is interested not only in the reduced C*-algebra of the inverse semigroup, but also in its boundary quotient. Let us describe the partial system corresponding to this quotient. Given a semilattice $E$, let $\widehat{E}_{\max}$ be the subset of $\widehat{E}$ consisting of those $\chi \in \widehat{E}$ such that $\menge{e \in E}{\chi(e) = 1}$ is maximal. We then set $\partial \widehat{E} \defeq \overline{\widehat{E}_{\max}} \subseteq \widehat{E}$. Now let $E$ be the semilattice of idempotents in an inverse semigroup $S$. As $\partial \widehat{E} \subseteq \widehat{E}$ is closed, we obtain a short exact sequence $0 \to I \to C_0(\widehat{E}) \to C_0(\partial \widehat{E}) \to 0$. Viewing $I$ as a subset of $C^*_{\lambda}(S)$, we form the ideal $\spkl{I}$ of $C^*_{\lambda}(S)$ generated by $I$. The boundary quotient in Exel's sense (see \cite{Ex2,Ex3,Ex4,EGS}) is given by $\partial C^*_{\lambda}(S) \defeq C^*_{\lambda}(S) / \spkl{I}$.

\blemma
\label{bdE}
Let $S$ be an inverse semigroup as above, with an idempotent pure partial homomorphism $\sigma$ from $S$ to a group $G$, with semilattice of idempotents $\widehat{E}$ and partial system $G \curvearrowright \widehat{E}$. Then $\partial \widehat{E}$ is $G$-invariant.
\elemma
\setlength{\parindent}{0.5cm} \setlength{\parskip}{0cm}
\bproof
Let us first show that for every $g \in G$, $g.(U_{g^{-1}} \cap \widehat{E}_{\max}) \subseteq U_g \cap \widehat{E}_{\max}$. Take $\chi \in \widehat{E}_{\max}$ with $\chi(s^*s) = 1$ for some $s \in S$ with $\sigma(s) = g$. Then $g.\chi(e) = \chi(s^*es)$. Assume that $g.\chi \notin \widehat{E}_{\max}$. This means that there is $\psi \in \widehat{E}_{\max}$ such that $\psi(e) = 1$ for all $e \in E$ with $\chi(e) = 1$, and there exists $f \in E$ with $\psi(f) = 1$ but $\chi(s^*fs) = 0$. Then $\psi \in U_g$ since $\chi(ss^*) = 1$, which implies $\psi(ss^*) = 1$. Consider $g^{-1}.\psi$ given by $g^{-1}.\psi(e) = \psi(ses^*)$. Then for every $e \in E$, $\chi(e) = 1$ implies $\chi(s^*ses^*s) = 1$, hence $\chi(s^*(ses^*)s) = 1$, so that $g^{-1}.\psi(e) = \psi(ses^*) = 1$. But $\chi(s^*fs) = 0$ and $g^{-1}.\psi(s^*fs) = \psi(ss^*fss^*) = \psi(f) = 1$. This contradicts $\chi \in \widehat{E}_{\max}$. Hence $g.(U_{g^{-1}} \cap \widehat{E}_{\max}) \subseteq U_g \cap \widehat{E}_{\max}$. To see that $g.(U_{g^{-1}} \cap \partial \widehat{E}) \subseteq U_g \cap \partial \widehat{E}$, let $\chi \in U_{g^{-1}} \cap \partial \widehat{E}$ and choose a net $(\chi_i)_i$ in $\widehat{E}_{\max}$ with $\lim_i \chi_i = \chi$. As $U_{g^{-1}}$ is open, we may assume that all the $\chi_i$ lie in $U_{g^{-1}}$. Then $g.\chi_i \in \widehat{E}_{\max}$, and $\lim_i g.\chi_i = g.\chi$. This implies $g.\chi \in \partial \widehat{E}$.
\eproof
\setlength{\parindent}{0cm} \setlength{\parskip}{0.5cm}

\blemma
\label{bdQ}
Let $S$ be an inverse semigroup with an idempotent pure partial homomorphism $\sigma$ from $S$ to an exact group $G$. Then the identification $C^*_{\lambda}(S) \cong C_0(\widehat{E}) \rtimes_r G$ from Proposition~\ref{S-EG} identifies $\partial C^*_{\lambda}(S)$ with $C_0(\partial \widehat{E}) \rtimes_r G$.
\elemma
\setlength{\parindent}{0.5cm} \setlength{\parskip}{0cm}
\bproof
It is easy to see that under the identification $C^*_{\lambda}(S) \cong C_0(\widehat{E}) \rtimes_r G$ from Proposition~\ref{S-EG}, $\spkl{I}$ corresponds to $I \rtimes_r G$. Since $G$ is exact, we obtain a short exact sequence $0 \to I \rtimes_r G \to C^*_{\lambda}(S) \to \partial C^*_{\lambda}(S) \to 0$ by \cite[Theorem~22.9]{Ex5}.
\eproof
\setlength{\parindent}{0cm} \setlength{\parskip}{0.5cm}

\bcor
\label{amenable-nuclear}
Let $S$ be an inverse semigroup with an idempotent pure partial homomorphism $\sigma$ from $S$ to a group $G$. If $G$ is amenable, then both $C^*_{\lambda}(S)$ and $\partial C^*_{\lambda}(S)$ are nuclear. If $G$ is exact, then both $C^*_{\lambda}(S)$ and $\partial C^*_{\lambda}(S)$ are exact.
\ecor
\setlength{\parindent}{0.5cm} \setlength{\parskip}{0cm}
\bproof
By Proposition~\ref{S-EG} and Lemma~\ref{bdQ}, both $C^*_{\lambda}(S)$ and $\partial C^*_{\lambda}(S)$ can be described as reduced partial crossed products by $G$, hence as reduced C*-algebras of Fell bundles over $G$. Therefore, our claims follow from \cite[Theorem~20.7, Theorem~25.10 and Theorem~25.12]{Ex5}.
\eproof
\setlength{\parindent}{0cm} \setlength{\parskip}{0.5cm}

\subsection{Graphs}
\label{sec-graphs}

Let $\cE = (\cE^0,\cE^1,r,s)$ be a (countable) graph with vertices $\cE^0$, edges $\cE^1$ and range and source maps $r,s: \: \cE^1 \to \cE^0$. Let $\cE^*$ be the set of finite paths in $\cE$, and let $l(\mu)$ denote the length of a path $\mu \in \cE^*$. The graph inverse semigroup $S_{\cE}$ is given by $S_{\cE} = \menge{(\mu,\nu) \in \cE^* \times \cE^*}{s(\mu) = s(\nu)} \cup \gekl{0}$, where $(\mu,\nu)^* = (\nu,\mu)$ and
$$
(\mu,\nu)(\zeta,\eta) = 
\bfa
(\mu,\nu'\eta) & {\rm if} \ \nu = \zeta \nu',\\
(\mu\zeta',\eta) & {\rm if} \ \zeta = \nu \zeta',\\
0 & {\rm else}.
\efa
$$
The semilattice $E_{\cE}$ of idempotents of $S_{\cE}$ is given by $\gekl{(\mu,\mu) \in \cE^* \times \cE^*} \cup \gekl{0}$, hence can be identified canonically with $\cE^* \cup \gekl{0}$. Multiplication in $E_{\cE}$ is given by
$$
\mu \cdot \nu = 
\bfa
\mu & {\rm if} \ \mu = \nu \mu',\\
\nu & {\rm if} \ \nu = \mu \nu',\\
0 & {\rm else}.
\efa
$$
Note that we write $\mu \nu$ for concatenation of paths and $\mu \cdot \nu$ for the product in $E_{\cE}$.

It is easy to see that $C^*_{\lambda}(S_{\cE})$ is canonically isomorphic to the Toeplitz C*-algebra of $\cE$, and that $\partial C^*_{\lambda}(S_{\cE})$ is canonically isomorphic to the graph C*-algebra of $\cE$.

\bremark
\label{convention-graphs}
Note that we are using the convention that the partial isometry $s_{\mu}$ for $\mu \in \cE^1$ has source projection $e_{s(\mu)}$ corresponding to the source of $\mu$, and range projection dominated by $e_{r(\mu)}$, the projection corresponding to the range of $\mu$. This is the same convention as in \cite{Web}, but different from the one in \cite{BCW,HS}.
\eremark
Let us construct an idempotent pure partial homomorphism on $S_{\cE}$. Let $\Fz_{\cE^1}$ be the free group generated by $\cE^1$. We view $\mu \in \cE^*$ with $l(\mu) \geq 1$ as elements in $\Fz_{\cE^1}$ in a canonical way. Define $\sigma: \: S_{\cE}\reg \to \Fz_{\cE^1}$ by setting $\sigma(\mu,\nu) = \mu \nu^{-1}$ if $l(\mu), l(\nu) \geq 1$, $\sigma(\mu,\nu) = \mu$ if $l(\mu) \geq 1$ and $l(\nu) = 0$, $\sigma(\mu,\nu) = \nu^{-1}$ if $l(\mu) = 0$ and $l(\nu) \geq 1$, and $\sigma(\mu,\nu) = e$ if $l(\mu) = l(\nu) = 0$. Here $e$ is the identity of $\Fz_{\cE^1}$. It is easy to check that $\sigma$ is an idempotent pure partial homomorphism from $S_{\cE}$ to $\Fz_{\cE^1}$. Clearly, $\sigma$ is the universal one, in the sense of \cite{MS}.

\bremark
Let us show how a modification of $\sigma$ produces an easy argument for nuclearity of graph C*-algebras. We write $\Fz_{\cE^1}^+$ for the free semigroup generated by $\cE^1$. Let $\Fz_2^+$ and $\Fz_2$ be the free semigroup and the free group on two generators. Clearly, there is a semigroup embedding $\Fz_{\cE^1}^+ \into \Fz_2^+$ ($\cE^1$ is countable by assumption). Moreover, by \cite{Hoch}, we have an embedding $\Fz_2^+ \into \Fz_2 / \Fz_2''$, where $\Fz_2''$ is the second commutator subgroup of $\Fz_2$. Hence we obtain an embedding $\Fz_{\cE^1}^+ \into \Fz_2^+ \into \Fz_2 / \Fz_2''$. By universal property of $\Fz_{\cE^1}$, we obtain a homomorphism $\Fz_{\cE^1} \to \Fz_2 / \Fz_2''$, which has to factorize as $\Fz_{\cE^1} \to \Fz_{\cE^1} / \Fz_{\cE^1}'' \to \Fz_2 / \Fz_2''$, such that the diagram
\bgloz
  \xymatrix{
  \Fz_{\cE^1}^+ \ar@{^{(}->}[d] \ar@{^{(}->}[r] & \Fz_2^+ \ar@{^{(}->}[r] & \Fz_2 / \Fz_2''\\
  \Fz_{\cE^1} \ar[rr] & & \Fz_{\cE^1} / \Fz_{\cE^1}'' \ar[u]
  }
\egloz
commutes. Thus, the canonical homomorphism $\Fz_{\cE^1}^+ \to \Fz_{\cE^1} / \Fz_{\cE^1}''$ is injective.
\setlength{\parindent}{0.5cm} \setlength{\parskip}{0cm}

Now let $\sigma''$ be the composition $S_{\cE}\reg \overset{\sigma}{\lori} \Fz_{\cE^1} \to \Fz_{\cE^1} / \Fz_{\cE^1}''$. $\sigma''$ is again a partial homomorphism, and $\sigma''$ is idempotent pure because $\Fz_{\cE^1}^+ \to \Fz_{\cE^1} / \Fz_{\cE^1}''$ is injective. Moreover, $\Fz_{\cE^1} / \Fz_{\cE^1}''$ is solvable, in particular amenable. Hence Corollary~\ref{amenable-nuclear} implies that both the Toeplitz C*-algebra as well as the graph C*-algebra of $\cE$ are nuclear.
\setlength{\parindent}{0cm} \setlength{\parskip}{0cm}
\eremark
Let us now come back to $\sigma: \: S_{\cE}\reg \to \Fz_{\cE^1}$, and describe the corresponding partial actions $\Fz_{\cE^1} \curvearrowright \widehat{E_{\cE}}$ and $\Fz_{\cE^1} \curvearrowright \partial \cE$, where $\partial \cE \defeq \partial \widehat{E_{\cE}}$.
\setlength{\parindent}{0.5cm} \setlength{\parskip}{0cm}

We start with $\Fz_{\cE^1} \curvearrowright \widehat{E_{\cE}}$. For $g \in \Fz_{\cE^1}$, $U_{g^{-1}}$ is empty unless $g \in \img(\sigma)$, i.e., $g = \alpha \beta^{-1}$ for some paths $\alpha, \beta \in (\cE^* \setminus \cE^0) \cup \gekl{e}$, and if $\alpha, \beta$ both lie in $\cE^* \setminus \cE^0$, then $s(\alpha) = s(\beta)$. For such $g$, $U_{g^{-1}}$ consists of those $\chi \in \widehat{E_{\cE}}$ such that there exists $(\mu,\nu) \in S_{\cE}\reg$ with $\chi(\nu) = 1$ and $\sigma(\mu,\nu) = \alpha \beta^{-1}$. For such $\chi$ with $(\mu,\nu)$ as above, $(g.\chi)(\zeta) = \chi(\nu \zeta')$ if $\zeta = \mu \zeta'$ and $(g.\chi)(\zeta) = 0$ otherwise. By Proposition~\ref{S-EG}, the Toeplitz C*-algebra of $\cE$ is canonically isomorphic to $C^*_{\lambda}(S_{\cE}) \cong C^*_r(\Fz_{\cE^1} \ltimes \widehat{E_{\cE}}) \cong C_0(\widehat{E_{\cE}}) \rtimes_r \Fz_{\cE^1}$.

Following \cite{Web} (see also \cite{EL}), the partial system $\Fz_{\cE^1} \curvearrowright \partial \cE$ can be described explicitly as follows: As a set, we identify $\partial \cE$ with $\cE^{\infty} \cup \menge{\alpha \in \cE^*}{s(\alpha) \notin \cE^0_0}$, where $\cE^{\infty}$ is the set of infinite paths in $\cE$. To describe the topology, let $\mu \in \cE^*$ and set
$$
  \cZ(\mu) = \menge{\nu \in \cE^* \cup \cE^{\infty}}{\nu = \mu \nu' \ {\rm for} \ {\rm some} \ \nu' \in \cE^* \cup \cE^{\infty}}.
$$
Given a finite subset $F$ of $r^{-1}(s(\mu))$, let $\cZ(\mu \setminus F) = \cZ(\mu) \setminus \bigcup_{\nu \in F} \cZ(\mu \nu)$. Then
$$
  \partial \cE \cap \cZ(\mu \setminus F), \ \mu \in \cE^*, \, F \subseteq r^{-1}(s(\mu)) \, {\rm finite}
$$
is a basis for the topology of $\partial \cE$. The $\Fz_{\cE^1}$-action is just given by the restriction of the partial system $\Fz_{\cE^1} \curvearrowright \widehat{E_{\cE}}$, which has been described above. By Lemma~\ref{bdQ}, the graph C*-algebra of $\cE$ is canonically isomorphic to $\partial C^*_{\lambda}(S_{\cE}) \cong C^*_r(\Fz_{\cE^1} \ltimes \partial \cE) \cong C_0(\partial \cE) \rtimes_r \Fz_{\cE^1}$ (compare \cite[Theorem~20.9]{Ex2}).
\setlength{\parindent}{0cm} \setlength{\parskip}{0.5cm}

It is known (see \cite[Proposition~2.3]{BCW}, and also Lemma~\ref{graphs-GPD}) that $\Fz_{\cE^1} \ltimes \partial \cE$ is topologically principal, or equivalently (see Corollary~\ref{TF-->tp}), that $\Fz_{\cE^1} \curvearrowright \partial \cE$ is topologically free, if and only if $\cE$ satisfies condition (L). Recall that $\cE$ satisfies condition (L) if every loop has an entry, i.e., for every $\mu = \mu_1 \cdots \mu_n \in \cE^*$ with $\mu_i \in \cE^1$, $n \geq 1$ and $s(\mu_n) = r(\mu_1)$, there is $\nu \in \cE^1$ and $1 \leq i \leq n$ with $r(\nu) = r(\mu_i)$ and $\nu \neq \mu_i$. 

Let us now relate our transformation groupoid to the groupoid attached to topological Markov shifts or graphs (see \cite{MM} and \cite{BCW}), and our notion of continuous orbit equivalence for partial systems to continuous orbit equivalence for topological Markov shifts or graphs (see \cite{MM} and \cite{BCW}). Once these relations are established, we will see that Theorem~\ref{COE-gpd-C} applied to graphs generalizes \cite[Theorem~2.3]{MM} and gives an alternative interpretation for \cite[Theorem~5.1]{BCW}.

Let $\cE$ be a graph as above. We compare the transformation groupoid $\Fz_{\cE^1} \ltimes \partial \cE$ with the groupoid $\cG_{\cE}$ from \cite[\S~2.3]{BCW} (see \cite[\S~2.2]{MM} for the case of topological Markov shifts). The groupoid $\cG_{\cE}$ is given by $\cG_{\cE} = \menge{(\alpha,n,\beta) \in \partial \cE \times \Zz \times \partial \cE}{n = k - l \ for \ k, l \in \Zz_+ \ and \ \sigma^k(\alpha) = \sigma^l(\beta)}$. Here, we identify $\partial \cE$ with $\cE^{\infty} \cup \menge{\alpha \in \cE^*}{s(\alpha) \notin \cE^0_0}$ and define for $\mu = \mu_1 \mu_2 \mu_3 \dotso$ in $\cE^{\infty}$ with $l(\mu) \geq 2$ ($\mu_i \in \cE^1$) $\sigma(\mu) = \mu_2 \mu_3 \dotso$, and $\sigma(\mu) = s(\mu)$ if $\mu \in \cE^1$. For $\mu \in \cE^0$, $\sigma$ is not defined. An equation like $\sigma^k(\alpha) = \sigma^l(\beta)$ always implicitly means that $\sigma^k(\alpha)$ and $\sigma^l(\beta)$ are defined, i.e., $l(\alpha) \geq k$ and $l(\beta) \geq l$. Moreover, we write $\Zz_+ = \menge{z \in \Zz}{z \geq 0}$.
\blemma
\label{graphs-GPD}
$\Fz_{\cE^1} \ltimes \partial \cE$ and $\cG_{\cE}$ are isomorphic as topological groupoids.
\elemma
\setlength{\parindent}{0.5cm} \setlength{\parskip}{0cm}
\bproof
It is easy to see that $\cG_{\cE} \to \Fz_{\cE^1} \ltimes \partial \cE, \, (\lambda \nu, n, \mu \nu) \ma (\lambda \mu^{-1}, \mu \nu)$ and $\Fz_{\cE^1} \ltimes \partial \cE \to \cG_{\cE}, \, (\lambda \mu^{-1}, \mu \nu) \ma (\lambda \nu, l(\lambda) - l(\mu), \mu \nu)$ are mutually inverse (continuous) groupoid homomorphisms. (Note that these expressions for the maps only make sense for $l(\mu), l(\nu) \geq 1$. But there is an obvious way to define these maps for $l(\mu) = 0$ or $l(\nu) = 0$.)
\eproof
\setlength{\parindent}{0cm} \setlength{\parskip}{0cm}

Let us also compare the notion of (continuous) orbit equivalence for graphs introduced in \cite[Definition~3.1]{BCW} (see also \cite[\S~2.1]{MM} for the case of topological Markov shifts) with continuous orbit equivalence of the corresponding partial systems. Let $\cE$ and $\cF$ be two graphs, with $\sigma_{\cE}$ and $\sigma_{\cF}$ as above, and partial systems $\Fz_{\cE^1} \curvearrowright \partial \cE$ and $\Fz_{\cF^1} \curvearrowright \partial \cF$. $\cE$ and $\cF$ are (continuously) orbit equivalent in the sense of \cite[Definition~3.1]{BCW} if there exists a homeomorphism $\varphi: \: \partial \cE \to \partial \cF$  together with continuous maps $k, l: \: \partial \cE \to \Zz_+$ and $k', l': \: \partial \cF \to \Zz_+$ such that
\bgl
\label{sFsE}
  \sigma_{\cF}^{k(\zeta)}(\varphi(\sigma_{\cE}(\zeta))) = \sigma_{\cF}^{l(\zeta)}(\varphi(\zeta)) \ for \ all \ \zeta \in \partial \cE,
\egl
\bgl
\sigma_{\cE}^{k'(\eta)}(\varphi^{-1}(\sigma_{\cF}(\eta))) = \sigma_{\cE}^{l'(\eta)}(\varphi^{-1}(\eta)) \ for \ all \ \eta \in \partial \cF.
\egl

\blemma
\label{graphs-COE}
If $\Fz_{\cE^1} \curvearrowright \partial \cE$ and $\Fz_{\cF^1} \curvearrowright \partial \cF$ are continuously orbit equivalent, then $\cE$ and $\cF$ are orbit equivalent.
\setlength{\parindent}{0.5cm} \setlength{\parskip}{0cm}

If $\cE$ and $\cF$ satisfy condition (L), then the converse holds, i.e., if $\cE$ and $\cF$ are orbit equivalent, then $\Fz_{\cE^1} \curvearrowright \partial \cE$ and $\Fz_{\cF^1} \curvearrowright \partial \cF$ are continuously orbit equivalent.
\setlength{\parindent}{0cm} \setlength{\parskip}{0.5cm}
\elemma
\setlength{\parindent}{0.5cm} \setlength{\parskip}{0cm}
\bproof
Assume that $\Fz_{\cE^1} \curvearrowright \partial \cE$ and $\Fz_{\cF^1} \curvearrowright \partial \cF$ are continuously orbit equivalent via $\varphi: \: \partial \cE \cong \partial \cF$ and continuous maps $a$, $b$ as in Definition~\ref{COE}. As $\partial \cE = \rukl{\bigsqcup_{\zeta \in \cE^1} \zeta \partial \cE} \sqcup \cE^0 \setminus \cE^0_0$, it suffices to define $k$ and $l$ on $\bigsqcup_{\zeta \in \cE^1} \zeta \partial \cE$ (on the remaining part, just set $k$ and $l$ to be $0$). For $\zeta \in \cE^1$ and $\zeta \xi \in \zeta \partial \cE$, we know that $a(\zeta^{-1},\zeta \xi)$ has the reduced form $\lambda \kappa^{-1}$, so that we can define $k$ and $l$ on $\zeta \partial \cE$ by setting $k(\zeta \xi) \defeq l(\lambda)$ and $l(\zeta \xi) \defeq l(\kappa)$. These maps $k$ and $l$ are obviously locally constant, hence continuous. Moreover, we know that $\varphi(\zeta^{-1}.(\zeta \xi)) = (\lambda \kappa^{-1}).\varphi(\zeta \xi)$, which means that there exists $\omega \in \partial \cF$ with $\varphi(\zeta \xi) = \kappa \omega$ and $\varphi(\xi) = \varphi(\zeta^{-1}.(\zeta \xi)) = \lambda \omega$. Therefore, $\sigma_{\cF}^{l(\lambda)}(\varphi(\sigma_{\cE}(\zeta \xi))) = \sigma_{\cF}^{l(\lambda)}(\varphi(\xi)) = \sigma_{\cF}^{l(\lambda)}(\lambda \omega) = \omega = \sigma_{\cF}^{l(\kappa)}(\kappa \omega) = \sigma_{\cF}^{l(\kappa)}(\varphi(\zeta \xi))$. Thus, \eqref{sFsE} holds. $k'$ and $l'$ are defined in a similar way, using $b$.
\setlength{\parindent}{0cm} \setlength{\parskip}{0.5cm}

Now assume that $\cE$ and $\cF$ satisfy condition (L), and suppose conversely that $\cE$ and $\cF$ are orbit equivalent. Because of condition (L), our partial systems $\Fz_{\cE^1} \curvearrowright \partial \cE$ and $\Fz_{\cF^1} \curvearrowright \partial \cF$ are topologically free. Therefore, to prove that they are continuously orbit equivalent, all we have to show is that for every $g \in \Fz_{\cE^1}$ and $x \in U_{g^{-1}}$, there exists an open neighbourhood $U$ of $x$ and $h \in \Fz_{\cF^1}$ such that $\varphi(g.\bar{x}) = h.\varphi(\bar{x})$ for all $\bar{x} \in U$, and analogously for $\varphi^{-1}$. In our case, since $\Fz_{\cE^1}$ is generated by $\cE^1$, it suffices to consider $g \in (\cE^1)^{-1}$. Take $\zeta \in \cE^1$ and $x \in \zeta \partial \cE$. Choose $\lambda, \kappa \in \cF^*$ with $l(\lambda) \geq k(x)$, $l(\kappa) \geq l(x)$, such that $\varphi(x) \in \kappa \partial \cF$ and $\varphi(\zeta^{-1}.x) \in \lambda \partial \cF$, and for all $\zeta \xi \in \zeta \partial \cE$ with $\varphi(\zeta \xi) \in \kappa \partial \cF$ and $\varphi(\xi) = \varphi(\zeta^{-1}.(\zeta \xi)) \in \lambda \partial \cF$, we have $k(\zeta \xi) = k(x)$ and $l(\zeta \xi) = l(x)$. Such $\lambda$ and $\kappa$ exist because $k$ and $l$ are continuous, hence locally constant. Set $U \defeq \menge{\zeta \xi \in \zeta \partial \cE}{\varphi(\zeta \xi) \in \kappa \partial \cF \ and \ \varphi(\xi) \in \lambda \partial \cF}$. $U$ is obviously an open neighbourhood of $x$. Set $k \defeq k(x)$ and $l \defeq l(x)$. For all $\zeta \xi \in U$, we have $\sigma_{\cF}^k(\varphi(\xi)) = \sigma_{\cF}^k(\varphi(\sigma_{\cE}(\zeta \xi))) = \sigma_{\cF}^l(\varphi(\zeta \xi))$. Thus, when we write $\lambda = \lambda' \lambda''$ with $l(\lambda') = k$ and $\kappa = \kappa' \kappa''$ with $l(\kappa') = l$, we even know that $(\lambda')^{-1}.\varphi(\xi) = \sigma_{\cF}^k(\varphi(\xi)) = \sigma_{\cF}^l(\varphi(\zeta \xi)) = (\kappa')^{-1}.\varphi(\zeta \xi)$, and hence $\varphi(\zeta^{-1}.(\zeta \xi)) = (\lambda' (\kappa')^{-1}).\varphi(\zeta \xi)$ for all $\zeta \xi \in U$. Thus $\varphi$ has the desired property. The proof for $\varphi^{-1}$ is analogous.
\eproof
\setlength{\parindent}{0cm} \setlength{\parskip}{0cm}

Lemma~\ref{graphs-GPD}, Lemma~\ref{graphs-COE} and Theorem~\ref{COE-gpd-C} imply the following
\setlength{\parindent}{0.5cm} \setlength{\parskip}{0cm}
\btheo
Let $\cE$ and $\cF$ be graphs. Consider the statements
\begin{enumerate}
\item[a)] $\Fz_{\cE^1} \ltimes \partial \cE$ and $\Fz_{\cF^1} \ltimes \partial \cF$ are isomorphic as topological groupoids,
\item[b)] $\cG_{\cE}$ and $\cG_{\cF}$ are isomorphic as topological groupoids,
\item[c)] there exists an isomorphism $\Phi: \: C^*(\cE) \overset{\cong}{\lori} C^*(\cF)$ with $\Phi(C^*(E_{\cE})) = C^*(E_{\cF})$,
\item[d)] $\Fz_{\cE^1} \curvearrowright \partial \cE$ and $\Fz_{\cF^1} \curvearrowright \partial \cF$ are continuously orbit equivalent,
\item[e)] $\cE$ and $\cF$ are orbit equivalent.
\end{enumerate}
We always have a) $\LRarr$ b), a) $\Rarr$ c) as well as a) $\Rarr$ d) $\Rarr$ e). If $\cE$ and $\cF$ satisfy condition (L), then all these statements are equivalent.
\etheo
\setlength{\parindent}{0cm} \setlength{\parskip}{0cm}
This generalizes \cite[Theorem~2.3]{MM} and gives an alternative proof for the corresponding parts of \cite[Theorem~5.1]{BCW}. Note that in \cite{BCW}, it is proven that a), b) and c) are always equivalent.
\setlength{\parindent}{0cm} \setlength{\parskip}{0.5cm}

\subsection{Subsemigroups of groups}

Let us discuss C*-algebras attached to subsemigroups of groups. Given a subsemigroup $P$ of a group $G$ (or any left cancellative semigroup), the left reduced semigroup C*-algebra $C^*_{\lambda}(P)$ of $P$ is defined as the sub-C*-algebra of $\cL(\ell^2 P)$ generated by the left multiplication operators $V_p: \: \ell^2 P \to \ell^2 P, \, \delta_x \ma \delta_{px}$. Here $\delta_x$ is the delta function in $x \in P$, and these $\delta_x$ form the canonical orthonormal basis of $\ell^2 P$. The canonical commutative subalgebra $D_{\lambda}(P)$ is given by $D_{\lambda}(P) = C^*_{\lambda}(P) \cap \ell^{\infty}(P)$. Here is an alternative description of $D_{\lambda}(P)$: Let $I_V$ be the inverse semigroup of partial isometries on $\ell^2 P$ generated by $V_p$, $p \in P$, i.e., $I_V = \menge{V_{p_1}^* V_{q_1} \dotsm V_{p_n}^* V_{q_n}}{n \in \Nz, \, p_i, q_i \in P}$. We can define a partial homomorphism $\sigma: \: I_V \setminus \gekl{0} \to G, \, V_{p_1}^* V_{q_1} \dotsm V_{p_n}^* V_{q_n} \ma p_1^{-1} q_1 \dotsm p_n^{-1} q_n$. To see that $\sigma$ is well-defined, note that every $V \in I_V$ has the property that there exists $g \in G$ such that for every $x \in P$, either $V \delta_x = 0$ or $V \delta_x = \delta_{gx}$. And $\sigma$ is defined in such a way that $\sigma(V) = g$. Now the closed linear span $\clspan(\sigma^{-1}(e))$ of $\sigma^{-1}(e)$, where $e$ is the identity of $G$, coincides with $D_{\lambda}(P)$. All this is explained in \cite{Li1,Li2} (see \cite[Remark~3.12]{Li1} for the alternative description of $D_{\lambda}(P)$).

Let us now describe the canonical partial action $G \curvearrowright D_{\lambda}(P)$. We first describe the dual action $\alpha^*$. For $g \in G$, let $D_{g^{-1}} \defeq \clspan(\menge{V^*V}{V \in I_V \setminus \gekl{0}, \, \sigma(V) = g})$. $D_{g^{-1}}$ is an ideal of $D_{\lambda}(P)$, and it follows from the alternative description of $D_{\lambda}(P)$ that $D_e = D_{\lambda}(P)$. We then define $\alpha^*_g$ as $\alpha^*_g : \: D_{g^{-1}} \to D_g, \, V^*V \to VV^*$ for $V \in I_V \setminus \gekl{0}$ with $\sigma(V) = g$. This is well-defined: If we view $\ell^2 P$ as a subspace $\ell^2 G$ and let $\lambda$ be the left regular representation of $G$, then every $V \in I_V \setminus \gekl{0}$ with $\sigma(V) = g$ satisfies $V = \lambda_g V^*V$. Therefore, $VV^* = \lambda_g V^*V \lambda_g^*$. This shows that $\alpha^*_g$ is just conjugation with the unitary $\lambda_g$. This also explains why $\alpha^*_g$ is an isomorphism. Of course, we can also describe the dual action $\alpha$. Set $\Omega_P \defeq \Spec(D_{\lambda}(P))$ and for every $g \in G$, let $U_{g^{-1}} \defeq \widehat{D_{g^{-1}}}$. It is easy to see that $U_{g^{-1}} = \menge{\chi \in \Omega_P}{\chi(V^*V) = 1 \ for \ some \ V \in I_V \setminus \gekl{0} \ with \ \sigma(V) = g}$. We then define $\alpha_g$ by setting $\alpha_g(\chi) \defeq \chi \circ \alpha^*_{g^{-1}}$.

\bprop
\label{P-DG}
There is a canonical isomorphism $C^*_{\lambda}(P) \cong D_{\lambda}(P) \rtimes_r G$ determined by $V_p \ma V_p V_p^* u_p$. Here we form the partial crossed product for the partial action $G \curvearrowright D_{\lambda}(P)$ defined above, and $u_g$ denote the canonical partial isometries in $D_{\lambda}(P) \rtimes_r G$.
\eprop
\setlength{\parindent}{0.5cm} \setlength{\parskip}{0cm}
\bproof
Our strategy is to describe both $C^*_{\lambda}(P)$ and $D_{\lambda}(P) \rtimes_r G$ as reduced (cross sectional) algebras of Fell bundles, and then to identify the underlying Fell bundles.
\setlength{\parindent}{0cm} \setlength{\parskip}{0.25cm}

Let us start with $C^*_{\lambda}(P)$. We have already defined $I_V$ and $\sigma$. Now we set $B_g \defeq \clspan(\sigma^{-1}(g))$ for every $g \in G$. We want to see that $(B_g)_{g \in G}$ is a grading for $C^*_{\lambda}(P)$, in the sense of \cite[Definition~3.1]{Ex1}. Conditions (i) and (ii) are obviously satisfied. For (iii), we use the faithful conditional expectation $\cE: \: C^*_{\lambda}(P) \onto D_{\lambda}(P) = B_e$ from \cite[\S~3.1]{Li1}. Given a finite sum $x = \sum_g x_g \in C^*_{\lambda}(P)$ of elements $x_g \in B_g$ such that $x = 0$, we conclude that $0 = x^* x = \sum_{g,h} x_g^* x_h$, and hence $0 = \cE(x^* x) = \sum_g x_g^* x_g$ (here we used that $\cE \vert_{B_g} = 0$ if $g \neq e$). This implies that $x_g = 0$ for all $g$. Therefore, the subspaces $B_g$ are independent. It is clear that the linear span of all the $B_g$ is dense in $C^*_{\lambda}(P)$. This proves (iii). If we let $\cB$ be the Fell bundle given by $(B_g)_{g \in G}$, then \cite[Proposition~3.7]{Ex1} implies $C^*_{\lambda}(P) \cong C^*_r(\cB)$ because $\cE: \: C^*_{\lambda}(P) \onto D_{\lambda}(P) = B_e$ is a faithful conditional expectation satisfying $\cE \vert_{B_e} = \id_{B_e}$ and $\cE \vert_{B_g} = 0$ if $g \neq e$.

Let us also describe $D_{\lambda}(P) \rtimes_r G$ as a reduced algebra of a Fell bundle. We denote by $W_g$ the partial isometry in $D_{\lambda}(P) \rtimes_r G$ corresponding to $g \in G$, and we set $B'_g \defeq D_g W_g$. Recall that we defined $D_{g^{-1}} = \clspan(\menge{V^*V}{V \in I_V \setminus \gekl{0}, \, \sigma(V) = g})$ earlier on. It is easy to check that $(B'_g)_{g \in G}$ satisfy (i), (ii) and (iii) in \cite[Definition~3.1]{Ex1}. Moreover, $B'_e = D_e = D_{\lambda}(P)$, and it follows immediately from the construction of the reduced partial crossed product that there is a faithful conditional expectation $D_{\lambda}(P) \rtimes_r G \onto D_{\lambda}(P) = B'_e$ which is identity on $B'_e$ and $0$ on $B'_g$ for $g \neq e$. Hence if we let $\cB'$ be the Fell bundle given by $(B'_g)_{g \in G}$, then \cite[Proposition~3.7]{Ex1} implies $D_{\lambda}(P) \rtimes_r G \cong C^*_r(\cB')$.

To identify $C^*_{\lambda}(P)$ and $D_{\lambda}(P) \rtimes_r G$, it now remains to identify $\cB$ with $\cB'$. We claim that the map $\lspan(\menge{V}{\sigma(V) = g}) \to \lspan(\menge{VV^* W_g}{\sigma(V) = g}), \, \sum_i \alpha_i V_i \ma \sum_i \alpha_i V_i V_i^* W_g$ is well-defined and extends to an isometric isomorphism $B_g \to B'_g$, for all $g \in G$.

All we have to show is that our map is isometric. We have
$\norm{\sum_i \alpha_i V_i}^2 = \norm{\sum_{i,j} \alpha_i \overline{\alpha_j} V_i V_j^*}_{D_{\lambda}(P)}$
and
$\norm{\sum_i \alpha_i V_i V_i^* W_g}^2 = \norm{\sum_{i,j} \alpha_i \overline{\alpha_j} V_i V_i^* V_j V_j^*}_{D_{\lambda}(P)}$.
Since $V_i = V_i V_i^* \lambda_g$ and $V_j^* = \lambda_{g^{-1}} V_j V_j^*$, we have $V_i V_j^* = V_i V_i^* \lambda_g \lambda_{g^{-1}} V_j V_j^* = V_i V_i^* V_j V_j^*$. Hence, indeed, $\norm{\sum_i \alpha_i V_i}^2 = \norm{\sum_i \alpha_i V_i V_i^* W_g}^2$, and we are done.

All in all, we have proven that $C^*_{\lambda}(P) \cong C^*_r(\cB) \cong C^*_r(\cB') \cong D_{\lambda}(P) \rtimes_r G$.
\eproof
\setlength{\parindent}{0cm} \setlength{\parskip}{0.5cm}

\bremark
A straightforward computation shows that actually, $V_p V_p^* u_p = u_p$ for all $p \in P$. Thus the isomorphism in the Proposition~\ref{P-DG} is given by $V_p \ma u_p$ for all $p \in P$.
\eremark

Our next goal is to write $C^*_{\lambda}(P)$ as a quotient of a C*-algebra of an inverse semigroup in a canonical way. Let $S \defeq I_l(P)$ be the inverse semigroup of partial bijections of $P$ generated by $p: \: P \overset{\cong}{\lori} P, \, x \ma px$ ($p \in P$). The semilattice of idempotents $E$ of $S$ is given by the set of constructible ideals $\cJ = \menge{p_1^{-1} q_1 \cdots p_n^{-1} q_n P}{p_i, q_i \in P} \cup \gekl{\emptyset}$. It is easy to see that $S$ is canonically isomorphic to the inverse semigroup $I_V$ constructed above. The homomorphism $I_V \setminus \gekl{0} \to G$ yields an idempotent pure partial homomorphism from $S$ to $G$ such that for every $s \in S\reg$, $s(x) = \sigma(s) \cdot x$ if $x \in \dom(s)$.
\setlength{\parindent}{0cm} \setlength{\parskip}{0.25cm}

As explained in \cite[Corollary~3.2.13]{Nor}, the isometry $\ell^2 P \to \ell^2 S\reg, \, \delta_p \ma \delta_p$ induces surjective homomorphisms $C^*(E) \onto D_{\lambda}(P)$ and $C^*_{\lambda}(S) \onto C^*_{\lambda}(P)$. The first surjection allows us to view $\Omega_P = \Spec(D_{\lambda}(P))$ as a closed subspace of $\widehat{E}$. More precisely, $\chi \in \widehat{E}$ lies in $\Omega_P$ if for all constructible ideals $X, X_1, \dotsc, X_n$ of $P$ with $X = \bigcup_{i=1}^n X_i$, $\chi(X) = 1$ implies $\chi(X_i) = 1$ for some $1 \leq i \leq n$.

The following lemma is easy to check:
\blemma
\label{P-S}
$\Omega_P$ is a $G$-invariant subspace with respect to the canonical $G$-action on $\widehat{E}$. The corresponding partial system $G \curvearrowright \Omega_P$ coincides with $\alpha$, the dual action of $\alpha^*: \: G \curvearrowright D_{\lambda}(P)$ constructed before Proposition~\ref{P-DG}.
\setlength{\parindent}{0.5cm} \setlength{\parskip}{0cm}

Moreover, the surjection $C^*_{\lambda}(S) \onto C^*_{\lambda}(P)$ from \cite[Corollary~3.2.13]{Nor} corresponds to the obvious map $C(\widehat{E}) \rtimes_r G \onto C(\Omega_P) \rtimes_r G$ under the identifications given by Proposition~\ref{S-EG} and Proposition~\ref{P-DG}.
\setlength{\parindent}{0cm} \setlength{\parskip}{0.25cm}
\elemma

Here is a sufficient condition for topological freeness of $G \curvearrowright \Omega_P$:
\bprop
If $P$ contains the identity $e \in G$, and if the group of units $P^*$ in $P$ is trivial, i.e., $P^* = \gekl{e}$, then $G \curvearrowright \Omega_P$ is topologically free.
\eprop
\setlength{\parindent}{0.5cm} \setlength{\parskip}{0cm}
\bproof
For $p \in P$, let $\chi_p \in \widehat{E}$ be defined by $\chi_p(X) = 1$ if and only if $p \in X$, for $X \in \cJ$. Obviously, $\chi_p \in \Omega_P$. It turns out that $\menge{\chi_p}{p \in P}$ is dense in $\Omega_P$. Basic open sets in $\Omega_P$ are of the form $U(X;X_1, \dotsc, X_n) = \menge{\chi \in \Omega_P}{\chi(X) = 1, \, \chi(X_i) = 0 \ {\rm for} \ {\rm all} \ 1 \leq i \leq n}$. Here $X, X_1, \dotsc, X_n$ are constructible ideals of $P$. Clearly, $U(X;X_1, \dotsc, X_n)$ is empty if $X = \bigcup_{i=1}^n X_i$. Thus, for a non-empty basic open set $U(X;X_1, \dotsc, X_n)$, we may choose $p \in X$ such that $p \notin \bigcup_{i=1}^n X_i$, and then $\chi_p \in U(X;X_1, \dotsc, X_n)$.
\setlength{\parindent}{0cm} \setlength{\parskip}{0.25cm}

Let $p \in P$ and $g \in G$ satisfy $g.\chi_p = \chi_p$. This equality only makes sense if $\chi_p \in U_{g^{-1}}$, i.e., there exists $Y \in \cJ$ with $Y \subseteq g^{-1} \cdot P$ and $p \in Y$. Then $g.\chi_p(X) = \chi_p(g^{-1} \cdot (X \cap g \cdot Y))$. So $g.\chi_p(X) = 1$ if and only if $p \in g^{-1} \cdot (X \cap g \cdot Y)$ if and only if $gp \in X \cap g \cdot Y$ if and only if $gp \in X$, while $\chi_p(X) = 1$ if and only if $p\in P$. Hence $g.\chi_p = \chi_p$ implies that for every $X \in \cJ$, $gp \in X$ if and only if $p \in X$. For $X = pP$, we obtain $gp \in pP$, and for $X = gpP$, we get $p \in gpP$. Hence there exist $x,y \in P$ with $gp = px$ and $p = gpy$. So $p = gpy = pxy$ and $gp = px = gpyx$. Thus $xy = yx = e$. Hence $x,y \in P^*$. Since $P^* = \gekl{e}$ by assumption, we must have $x = y = e$, and hence $gp = p$. This implies $g = e$. In other words, for every $e \neq g \in G$, $g.\chi_p \neq \chi_p$ for all $p \in P$ such that $\chi_p \in U_{g^{-1}}$. Hence it follows that $\menge{\chi \in U_{g^{-1}}}{g.\chi \neq \chi}$ contains $\menge{\chi_p \in U_{g^{-1}}}{p \in P}$, and the latter set is dense in $U_{g^{-1}}$ as $\menge{\chi_p}{p \in P}$ is dense in $\Omega_P$.
\eproof
\setlength{\parindent}{0cm} \setlength{\parskip}{0.5cm}
Coming back to the comparison of $C^*_{\lambda}(S)$ and $C^*_{\lambda}(P)$, it was shown in \cite[Theorem~3.2.14]{Nor} that the following are equivalent:
\setlength{\parindent}{0.5cm} \setlength{\parskip}{0cm}
\begin{itemize}
\item the canonical map $C^*_{\lambda}(S) \onto C^*_{\lambda}(P)$ is injective,
\item the canonical map $C^*(E) \onto D_{\lambda}(P)$ is injective,
\item $\cJ$ is independent.
\end{itemize}
Recall that $\cJ$ is called independent if for all $X, X_1, \dotsc, X_n \in \cJ$, $X = \bigcup_{i=1}^n X_i$ implies $X = X_i$ for some $1 \leq i \leq n$.
\setlength{\parindent}{0cm} \setlength{\parskip}{0.5cm}

In view of \cite{Li1,Li2}, it makes sense to view the full C*-algebra $C^*(S)$ of the inverse semigroup $S = I_l(P)$ as the full semigroup C*-algebra of $P$. Thus we set $C^*(P) \defeq C^*(S)$, and let $\lambda: \: C^*(P) \onto C^*_{\lambda}(P)$ be the composite $C^*(P) = C^*(S) \onto C^*_{\lambda}(S) \onto C^*_{\lambda}(P)$. Because of our descriptions as transformation groupoids (see Proposition~\ref{S-EG}), the following generalizations of \cite[Theorem~6.1]{Li2} are immediate.
\setlength{\parindent}{0.5cm} \setlength{\parskip}{0cm}
\btheo
\label{gen-Li1}
If $G \ltimes \widehat{E}$ or $G \ltimes \Omega_P$ is amenable, then $\lambda: \: C^*(P) \onto C^*_{\lambda}(P)$ is an isomorphism if and only if $P$ is independent.
\etheo
\btheo
\label{gen-Li2}
Consider the following statements:
\begin{enumerate}
\item[(i)] $C^*(P)$ is nuclear,
\item[(ii)] $C^*_{\lambda}(P)$ is nuclear,
\item[(iii)] $G \ltimes \Omega_P$ is amenable,
\item[(iv)] $\lambda: \: C^*(P) \onto C^*_{\lambda}(P)$ is an isomorphism.
\end{enumerate}
We always have (i) $\Rarr$ (ii) $\LRarr$ (iii). If $P$ is independent, then (iii) $\Rarr$ (iv) and (iii) $\Rarr$ (i), so that (i), (ii) and (iii) are equivalent.
\etheo
\setlength{\parindent}{0cm} \setlength{\parskip}{0.5cm}
\bcor
If $P$ is a subsemigroup of an amenable group $G$, then (i), (ii) and (iii) from Theorem~\ref{gen-Li2} hold, and (iv) is true if and only if $P$ is independent.
\ecor
\setlength{\parindent}{0.5cm} \setlength{\parskip}{0cm}
\bproof
By Proposition~\ref{S-EG}, both $C^*(P)$ and $C^*_{\lambda}(P)$ can be written as (full or reduced) partial crossed products by $G$, and hence as (full or reduced) C*-algebras of Fell bundles over $G$. Therefore, by \cite[Theorem~20.7 and Theorem~25.10]{Ex5}, $C^*(P)$ and $C^*_{\lambda}(P)$ are nuclear. (iii) holds because we know that (ii) $\LRarr$ (iii). Finally, our claim about (iv) follows from Theorem~\ref{gen-Li1}. 
\eproof
\setlength{\parindent}{0cm} \setlength{\parskip}{0cm}
Note that the statement in the last corollary was also obtained in \cite{Nor}.
\setlength{\parindent}{0cm} \setlength{\parskip}{0.25cm}

Now set $\partial \Omega_P \defeq \partial \widehat{E}$, where $E$ is the semilattice of idempotents of $S = I_l(P)$.
\blemma
\label{bdOmegaInOmega}
We have $\partial \Omega_P \subseteq \Omega_P$.
\elemma
\setlength{\parindent}{0.5cm} \setlength{\parskip}{0cm}
\bproof
Let $X, X_1, \dotsc, X_n \in \cJ$ satisfy $X = \bigcup_{i=1}^n X_i$. Then for $\chi \in \widehat{E}_{\max}$, $\chi(X_i) = 0$ implies that there exists $X_i' \in \cJ$ with $\chi(X_i') = 1$ and $X_i \cap X_i' = \emptyset$. Thus if $\chi(X_i) = 0$ for all $1 \leq i \leq n$, then let $X_i'$, $1 \leq i \leq n$ be as above. Then for $X' = \bigcap_{i=1}^n X_i'$, $\chi(X') = 1$ and $X \cap X' = \emptyset$. Thus $\chi(X) = 0$. This shows $\widehat{E}_{\max} \subseteq \Omega_P$. As $\Omega_P$ is closed, we conclude that $\partial \widehat{E} \subseteq \Omega_P$.
\eproof
\setlength{\parindent}{0cm} \setlength{\parskip}{0.25cm}

As $\partial \Omega_P$ is $G$-invariant, the following makes sense.
\bdefin
The boundary quotient of $C^*_{\lambda}(P)$ is given by $\partial C^*_{\lambda}(P) \defeq C(\partial \Omega_P) \rtimes_r G$.
\edefin
By construction and because of Proposition~\ref{P-DG}, there is a canonical projection $C^*_{\lambda}(P) \onto \partial C^*_{\lambda}(P)$. Morever, if $P$ is a subsemigroup of an exact group $G$, then $\partial C^*_{\lambda}(P)$ is isomorphic to $\partial C^*_{\lambda}(S) \cong C(\partial \Omega_P) \rtimes_r G$. Let us now generalize the results in \cite[\S~7.3]{Li2}.
\blemma
\label{Lem_bd-min}
$\partial \widehat{E} = \partial \Omega_P$ is the minimal non-empty closed $G$-invariant subspace of $\widehat{E}$.
\elemma
\setlength{\parindent}{0.5cm} \setlength{\parskip}{0cm}
\bproof
Let $C \subseteq \widehat{E}$ be non-empty, closed and $G$-invariant. Let $\chi \in \widehat{E}_{\max}$ be arbitrary, and choose $X \in \cJ$ with $\chi(X) = 1$. Choose $p \in X$ and $\chi \in C$. As $U_{p^{-1}} = \widehat{E}$, we can form $p.\chi$, and we know that $p.\chi \in C$. We have $p.\chi(pP) = \chi(P) = 1$, so that $p.\chi(X) = 1$ as $p \in X$ implies $pP \subseteq X$ ($X$ is a right ideal). Set $\chi_X \defeq p.\chi$. Consider the net $(\chi_X)_X$ indexed by $X \in \cJ$ with $\chi(X) = 1$, ordered by inclusion. Passing to a convergent subnet if necessary, we may assume that $\lim_X \chi_X$ exists. But it is clear because of $\chi \in \widehat{E}_{\max}$ that $\lim_X \chi_X = \chi$. As $\chi_X \in C$ for all $X$, we deduce that $\chi \in C$. Thus $\widehat{E}_{\max} \subseteq C$, and hence $\partial \widehat{E} \subseteq C$.
\eproof
\setlength{\parindent}{0cm} \setlength{\parskip}{0.25cm}

In particular, $\partial \Omega_P$ is the minimal non-empty closed $G$-invariant subspace of $\Omega_P$. Another immediate consequence is
\bcor
The transformation groupoid $G \ltimes \partial \Omega_P$ is minimal.
\ecor

To discuss topological freeness of $G \curvearrowright \partial \Omega_P$, let
$$G_0 = \menge{g \in G}{X \cap g \cdot P \neq \emptyset \neq X \cap g^{-1} \cdot P \ {\rm for} \ {\rm all} \ \emptyset \neq X \in \cJ},$$
as in \cite[\S~7.3]{Li2}. Clearly, $G_0 = \menge{g \in G}{pP \cap g \cdot P \neq \emptyset \neq pP \cap g^{-1} \cdot P \ {\rm for} \ {\rm all} \ p \in P}$. Furthermore, \cite[Lemma~7.19]{Li2} shows that $G_0$ is a subgroup of $G$.
\bprop
\label{Prop_bd-tf}
$G \curvearrowright \partial \Omega_P$ is topologically free if and only if $G_0 \curvearrowright \partial \Omega_P$ is topologically free. 
\eprop
\setlength{\parindent}{0.5cm} \setlength{\parskip}{0cm}
\bproof
\an{$\Rarr$} is clear. For \an{$\Larr$}, assume that $G_0 \curvearrowright \partial \Omega_P$ is topologically free, and suppose that $G \curvearrowright \partial \Omega_P$ is not topologically free, i.e., there exists $g \in G$ and $U \subseteq U_{g^{-1}} \cap \partial \Omega_P$ such that $g.\chi = \chi$ for all $\chi \in U$. As $\overline{\widehat{E}_{\max}} = \partial \Omega_P$, we can find $\chi \in U_{g^{-1}} \cap \widehat{E}_{\max}$ with $g.\chi = \chi$.
\setlength{\parindent}{0cm} \setlength{\parskip}{0.25cm}

For every $X \in \cJ$ with $\chi(X) = 1$, choose $x \in X$ and $\psi_X \in \widehat{E}_{\max}$ with $\psi_X(xP) = 1$, so that $\psi_X(X) = 1$. Consider the net $(\psi_X)_X$ indexed by $X \in \cJ$ with $\chi(X) = 1$, ordered by inclusion. Passing to a convergent subnet if necessary, we may assume that $\lim_X \psi_X = \chi$. As $U$ is open, we may assume that $\psi_X \in U$ for all $X$. Then $\psi_X(xP) = 1$ implies that $\psi_X \in U_x \cap U$.

Hence for sufficiently small $X \in \cJ$ with $\chi(X) = 1$, there exists $x \in X$ such that $x^{-1}.(U_x \cap U)$ is a non-empty open subset of $\partial \Omega_P$. We conclude that $(x^{-1}gx).\psi = \psi$ for all $\psi \in x^{-1}.(U_x \cap U)$. This implies that $x^{-1}gx \notin G_0$ as $G_0 \curvearrowright \partial \Omega_P$ is topologically free. So there exists $p \in P$ with $pP \cap x^{-1}gx \cdot P = \emptyset$ or $pP \cap x^{-1}g^{-1}x \cdot P = \emptyset$. Let $\chi_X \in \widehat{E}_{\max}$ satisfy $\chi_X(xpP) = 1$. If $pP \cap x^{-1}gx \cdot P = \emptyset$, then $xpP \cap gx \cdot P = \emptyset$, so that $xpP \cap g^{-1}xp \cdot P = \emptyset$. Hence $g.\chi_X \neq \chi_X$ if $\chi_X \in U_{g^{-1}}$. If $pP \cap x^{-1}g^{-1}x \cdot P = \emptyset$, then $xpP \cap g^{-1}x \cdot P = \emptyset$, so that $xpP \cap g^{-1}xp \cdot P = \emptyset$. Again, $g.\chi_X \neq \chi_X$ if $\chi_X \in U_{g^{-1}}$.

For every sufficiently small $X \in \cJ$ with $\chi(X) = 1$, we can find $x \in X$ and $\chi_X$ as above. Hence we can consider the net $(\chi_X)_X$ as above, and assume after passing to a convergent subnet that $\lim_X \chi_X = \chi$. As $\chi \in U \subseteq U_{g^{-1}} \cap \partial \Omega_P$, it follows that $\chi_X \in U \subseteq U_{g^{-1}} \cap \partial \Omega_P$ for sufficiently small $X$. So we obtain $g.\chi_X \neq \chi_X$, although $g$ acts trivially on $U$. This is a contradiction.
\eproof

\bcor
\label{Cor_bd-tf}
If $G_0 \curvearrowright \partial \Omega_P$ is topologically free, then $\partial C^*_{\lambda}(P)$ is simple.
\ecor
\bproof
This follows from Lemma~\ref{Lem_bd-min}, Proposition~\ref{Prop_bd-tf} and \cite[Chapter~II, Proposition~4.6]{R}. 
\eproof

\btheo
\label{THM_sgp-rep}
Let $P = \spkl{\cS, \cR}^+$ be a monoid given by a positive r-complete presentation $(\cS,\cR)$ in the sense of \cite{Deh}. Assume that for all $u \in \cS$, there is $v \in \cS$ such that $\cR$ does not contain any relation of the form $u \cdots = v \cdots$. Also, suppose that $P$ embeds into a group $G$ such that $(G,P)$ is quasi-lattice ordered in the sense of \cite{Ni}. Then $G_0 = \gekl{e}$ and $\partial C^*_{\lambda}(P)$ is simple.
\etheo
\bproof
In view of Corollary~\ref{Cor_bd-tf}, it suffices to prove $G_0 = \gekl{e}$. Let $g \in G_0$. Assume that $gP \cap P \neq P$. Then $g \in G_0$ implies that this intersection is not empty. Hence, we must have $gP \cap P = pP$ for some $p \in P$ because $(G,P)$ is quasi-lattice ordered. As $p \neq e$, there exists $u \in \Sigma$ with $pP \subseteq uP$. By assumption, there exists $v \in \Sigma$ such that no relation in $R$ is of the form $u \cdots = v \cdots$. By r-completeness, $uP \cap vP = \emptyset$ (see \cite[Proposition~3.3]{Deh}), so that $gP \cap vP = \emptyset$. This contradicts $g \in G_0$. Hence, we must have $gP \cap P = P$, and similarly, $g^{-1} P \cap P = P$. These two equalities imply $g \in P^*$. But $P^* = \gekl{e}$ because $(G,P)$ is quasi-lattice ordered. Thus $g = e$. 
\eproof

\bremark
By going over to the opposite semigroup, we obtain analogues results for the right versions $C^*_{\rho}(P)$ and $\partial C^*_{\rho}(P)$.
\eremark

\bexamples
Theorem~\ref{THM_sgp-rep} implies that for every graph-irreducible right-angled Artin monoid $A_{\Gamma}^+$ in the sense of \cite{CL}, $\partial C^*_{\lambda}(A_{\Gamma}^+)$ is simple. Also, for the Thompson monoid $F^+ = \spkl{x_0, x_1, \dotsc \ \vert \ x_n x_k = x_k x_{n+1} \ {\rm for} \ k<n}^+$, we get that $\partial C^*_{\rho}(F^+)$ is simple.

Corollary~\ref{Cor_bd-tf} and \cite[Corollary~5.10]{Li3} imply that for every countable Krull ring $R$ with $\cP(R)_{\tinf} \neq \emptyset$ or $\cP(R)_{\tfin}$ infinite (see \cite{Li3} for details), $\partial C^*_{\lambda}(R \rtimes R\reg)$ is simple. 
\eexamples
\setlength{\parindent}{0cm} \setlength{\parskip}{0.25cm}

\section{Purely infinite groupoids}
\label{pi}

Our goal is to exhibit examples of purely infinite groupoids. More precisely, we study groupoids attached to graphs, groupoids corresponding to boundary quotients of semigroup C*-algebras, and groupoids underlying semigroup C*-algebras of $ax+b$-semigroups.

The following lemma is easy to see, and will be used several times:
\blemma
\label{co=disjointunion}
Let $E$ be a semilattice. Every compact open subset $A \subseteq \widehat{E}$ can be written as a disjoint union $A = \bigsqcup_{i=1}^m U_i$ of basic open sets $U_i$ of the form
$$U(e;e_1, \dotsc, e_n) = \menge{\chi \in \widehat{E}}{\chi(e) = 1, \chi(e_1) = \dotso = \chi(e_n) = 0}, \ e, e_1, \dotsc, e_n \in E.$$
\elemma
This lemma will be helpful because given a partial system $G \curvearrowright \widehat{E}$ (or $G \curvearrowright X$ for any $G$-invariant subspace $X \subseteq \widehat{E}$), and we want to show that every compact open subset is $(G,\cC \cO)$-paradoxical in the sense of \cite[Definition~4.3]{GS}, then it suffices to show this for basic open sets of the form $U(e;e_1, \dotsc, e_n)$. Here $\cC \cO$ is the set of compact open subsets of $\widehat{E}$ or $X$.

\subsection{Graphs}

Let us first recall a necessary and sufficient condition from \cite{HS}, in terms of graphs, for pure infiniteness of graph C*-algebras. Let $\cE = (\cE^0,\cE^1,r,s)$ be a graph, and we use the same notation as in \S~\ref{sec-graphs}. Note that since our notation differs from the one in \cite{HS} (see Remark~\ref{convention-graphs}), we have to reverse all the arrows in the condition for pure infiniteness.

For $v, w \in \cE^0$, we write $w \leftarrow v$ if there exists a path $\mu \in \cE^*$ from $v$ to $w$, i.e., with  $r(\mu) = w$ and $s(\mu) = v$, and we write $w \not\leftarrow v$ if there is no such path. For $v \in \cE^0$, let $\Omega(v) = \menge{w \in \cE^0}{w \neq v, \, w \not\leftarrow v}$. $v \in \cE^0$ is called a breaking vertex if $\abs{r^{-1}(v)} = \infty$ and $0 < \abs{r^{-1}(v) \setminus s^{-1}(\Omega(v))} < \infty$.
\setlength{\parindent}{0.5cm} \setlength{\parskip}{0cm}

Moreover, we call $\mu \in \cE^*$ a loop if $l(\mu) \geq 1$ and $r(\mu) = s(\mu)$. We say that $\cE$ satisfies condition (K) if for all $v \in \cE^0$, whenever there exists a loop $\mu$ with $r(\mu) = v = s(\mu)$, there exists another loop $\mu'$ with $r(\mu') = v = s(\mu')$ and $\mu \cdot \mu' = 0$ in $E_{\cE}$.

Furthermore, a subset $M \subseteq \cE^0$ is called a maximal tail if
\begin{itemize}
\item for every $v \in \cE^0$, whenever there is $w \in M$ with $v \leftarrow w$, then $v \in M$;
\item for every $v \in M$ with $0 < \abs{r^{-1}(v)} < \infty$, there exists $\eta \in \cE^1$ with $r(\eta) = v$ and $s(\eta) \in M$;
\item for all $v,w \in M$, there exist $y \in M$ with $v \leftarrow y$ and $w \leftarrow y$.
\end{itemize}
We say that $v \in \cE^0$ connects to a loop if there is a loop $\mu$ with $r(\mu) = w = s(\mu)$ with $v \leftarrow w$.

\cite[Theorem~2.3]{HS} says that $C^*(\cE)$ is purely infinite if and only if $\cE$ satisfies the following condition:
\begin{itemize}
\item[(PI)] There exist no breaking vertices in $\cE$, $\cE$ satisfies condition (K), and every vertex in each maximal tail $M$ connects to a loop in $M$.
\end{itemize}
\setlength{\parindent}{0cm} \setlength{\parskip}{0.5cm}

\blemma
\label{inf-loops}
Suppose that $\cE$ satisfies condition (PI). Let $v \in \cE^0$ satisfy $\abs{r^{-1}(v)} = \infty$. Then there exist infinitely many loops $\mu_1, \mu_2, \dotsc \in \cE^*$ with $r(\mu_i) = v = s(\mu_i)$ of the form $\mu_i = \zeta_i \eta_i$ for pairwise distinct $\zeta_i \in \cE^1$.
\elemma
\setlength{\parindent}{0.5cm} \setlength{\parskip}{0cm}
\bproof
Since $v$ is not a breaking vertex, we must have $\abs{r^{-1}(v) \setminus s^{-1}(\Omega(v))} \in \gekl{0,\infty}$.

If $\abs{r^{-1}(v) \setminus s^{-1}(\Omega(v))} = \infty$, then there exist infinitely many pairwise distinct $\zeta_i \in \cE^1$ with $r(\zeta_i) = v$, and there is $\eta_i \in \cE^*$ with $r(\eta_i) = s(\zeta_i)$, $s(\eta_i) = v$. Then $\zeta_i \eta_i$ are the required loops.

If $\abs{r^{-1}(v) \setminus s^{-1}(\Omega(v))} = 0$, then consider the subset $M \defeq \menge{w \in \cE^0}{w \leftarrow v}$. It is easy to see that $M$ is a maximal tail, with $v \in M$. By condition (PI), $v$ has to connect to a loop $\mu$ in $M$. This means that there is $w \in \cE^0$ such that $w$ lies on $\mu$ and $v \leftarrow w$. But then, $w$ has to lie in $M$, so that $w \leftarrow v$ by definition of $M$; hence $v$ lies on a loop $v \leftarrow w \leftarrow v$. Let $\zeta_1 \cdots \zeta_n$ be such a loop. Then $\zeta_1 \in r^{-1}(v)$, but $\zeta_1 \notin s^{-1}(\Omega(v))$ as $s(\zeta_1) \leftarrow v$, for instance via $\zeta_2 \cdots \zeta_n$. Hence $0 < \abs{r^{-1}(v) \setminus s^{-1}(\Omega(v))}$, which is a contradiction.
\eproof
\setlength{\parindent}{0cm} \setlength{\parskip}{0cm}

Let $\Fz \ltimes \partial \cE$ be the transformation groupoid attached to $\Fz \curvearrowright \partial \cE$ as in \S~\ref{sec-graphs}.
\btheo
\label{PI-->pi}
If $\cE$ satisfies condition (PI), then $\Fz \ltimes \partial \cE$ is purely infinite.
\etheo
\setlength{\parindent}{0.5cm} \setlength{\parskip}{0cm}
\bproof
By Lemma~\ref{co=disjointunion}, it suffices to show that basic open sets of the form $U(\mu;\mu_1, \dotsc, \mu_n)$ are $(\Fz,\cC \cO)$-paradoxical. By passing even further to finite unions if necessary, it suffices to treat the basic open sets $U(\mu;\mu_1, \dotsc, \mu_n)$ with $\abs{r^{-1}(s(\mu))} = \infty$ (and $\mu_1, \dotsc, \mu_n \leq \mu$), or basic open sets of the form $U(\mu)$.

Let us consider the first case. Set $U \defeq U(\mu;\mu_1, \dotsc, \mu_n)$. Let $v = s(\mu)$. Since $\abs{r^{-1}(v)} = \infty$, Lemma~\ref{inf-loops} tells us that there are infinitely many loops $\mu_1, \mu_2, \dotsc \in \cE^*$ with $r(\mu_i) = v = s(\mu_i)$ of the form $\mu_i = \zeta_i \eta_i$ for pairwise distinct $\zeta_i \in \cE^1$. Hence we can find, among the $\zeta_i$, edges $\zeta, \zeta' \in \cE^1$ with $r(\zeta) = r(\zeta') = v$ and $\mu \zeta \cdot \mu_1 = \dotso = \mu \zeta \cdot \mu_n = 0 = \mu \zeta' \cdot \mu_1 = \dotso = \mu \zeta' \cdot \mu_n = 0$. Therefore, for $\chi \in \partial \cE$, $\chi(\mu \zeta) = 1$ implies $\chi(\mu_i) = 0$ for all $1 \leq i \leq n$, so that $\chi \in U$, and similarly for $\zeta'$. Moreover, $\mu \zeta \cdot \mu \zeta' = 0$ yields that $\chi(\mu \zeta) = 1$ implies $\chi(\mu \zeta') = 0$ and vice versa. Setting $g = \mu \zeta \mu^{-1} \in \Fz$, $h = \mu \zeta' \mu^{-1} \in \Fz$, we obtain $g.U \subseteq U$, $h.U \subseteq U$ and $g.U \cap h.U = \emptyset$.

In the second case, let $v = s(\mu)$. By Lemma~\ref{inf-loops} and condition (PI), we can find $\alpha = \alpha_1 \cdots \alpha_n \in \cE^*$ with $v = r(\alpha)$ and $\abs{r^{-1}(\alpha_i)} < \infty$ ($1 \leq i \leq n$) such that there is a loop $\beta \in \cE^*$ with $r(\beta) = s(\alpha) = s(\beta)$. Therefore, by passing to finite unions if necessary, we may assume by condition (K) that there exist two loops $\zeta$ and $\zeta'$ with $r(\zeta) = v = s(\zeta)$, $r(\zeta') = v = s(\zeta')$ and $\zeta \cdot \zeta' = 0$. Again, setting $g = \mu \zeta \mu^{-1} \in \Fz$, $h = \mu \zeta' \mu^{-1} \in \Fz$, we obtain $g.U \subseteq U$, $h.U \subseteq U$ and $g.U \cap h.U = \emptyset$.
\eproof
\setlength{\parindent}{0cm} \setlength{\parskip}{0cm}
It is clear that condition (PI) implies that $\Fz \curvearrowright \partial \cE$ is residually topologically free, in the sense of \cite[Definition~3.4~(ii)]{GS}. Moreover, \cite[Theorem~4.4]{GS} tells us that if $\Fz \curvearrowright \partial \cE$ is residually topologically free and if every compact open subset of $\partial \cE$ is $(\Fz,\cC \cO)$-paradoxical, then $C^*(\cE) \cong C(\partial \cE) \rtimes_r \Fz$ must be purely infinite (using that $\Fz$ is free and \cite[Theorem~22.9]{Ex5}). All in all, we obtain from Theorem~\ref{PI-->pi}:
\bcor
\label{graph_pi}
For a graph $\cE$, $\Fz \curvearrowright \partial \cE$ is residually topologically free and purely infinite if and only if $C^*(\cE)$ is purely infinite.
\ecor

\subsection{Boundary quotients of semigroup C*-algebras}

Let $P$ be a subsemigroup of a group $G$. Let $E$ be the semilattice of idempotents of $S = I_l(P)$. We write $\Omega$ for $\Omega_P$, $\Omega_{\max}$ for $\widehat{E}_{\max}$ and $\partial \Omega$ for $\partial \Omega_P = \partial \widehat{E}$.
\btheo
If there exist $p, q \in P$ with $pP \cap qP = \emptyset$, then $G \ltimes \partial \Omega$ is purely infinite.
\etheo
\setlength{\parindent}{0.5cm} \setlength{\parskip}{0cm}
\bproof
Let $U = \menge{\psi \in \partial \Omega}{\psi(X) = 1, \psi(X_1) = \dotso = \psi(X_n) = 0}$ be a basic open subset for some $X, X_1, \dotsc, X_n \in \cJ$. By Lemma~\ref{co=disjointunion}, it suffices to show that $U$ is $(G,\cC \cO)$-paradoxical. Since $\Omega_{\max}$ is dense in $\partial \Omega$, there exists $\chi \in \Omega_{\max}$ with $\chi \in U$. As $\chi$ lies in $\Omega_{\max}$, $\chi(X_i) = 0$ implies that there exists $Y_i \in \cJ$ with $X_i \cap Y_i = \emptyset$ and $\chi(Y_i) = 1$. Let $Y \defeq X \cap \bigcap_{i=1}^n Y_i$. Certainly, $Y \neq \emptyset$ as $\chi(Y) = 1$. Moreover, for every $\psi \in \partial \Omega$, $\psi(Y) = 1$ implies $\psi \in U$. Now choose $x \in Y$. By assumption, we can find $p, q \in P$ with $pP \cap qP = \emptyset$. For $\psi \in \partial \Omega$, $xp.\psi(xpP) = \psi(P) = 1$. Similarly, for all $\psi \in \partial \Omega$, $xq.\psi(xqP) = 1$. Thus $xp.U \subseteq xp.\partial \Omega \subseteq U$, $xq.U \subseteq xq.\partial \Omega \subseteq U$ and $(xp.U) \cap (xq.U) \subseteq (xp.\partial \Omega) \cap (xq.\partial \Omega) = \emptyset$ since $xpP \cap xqP = \emptyset$.
\eproof
\setlength{\parindent}{0cm} \setlength{\parskip}{0cm}
This strengthens and explains \cite[Corollary~7.23]{Li2}.
\setlength{\parindent}{0cm} \setlength{\parskip}{0.25cm}

\subsection{Groupoids from $ax+b$-semigroups}

Let $R$ be an integral domain with quotient field $K$. Consider the $ax+b$-semigroup $R \rtimes R\reg$, viewed as a subsemigroup of $K \rtimes K\reg$. Let $E$ be the semilattice of idempotents in $S = I_l(R \rtimes R\reg)$.
\btheo
Let $R$ be an integral domain with ${\rm Jac}(R) = (0)$ and $R\reg \neq R^*$. Then $(K \rtimes K\reg) \ltimes \widehat{E}$ is purely infinite.
\etheo
\setlength{\parindent}{0.5cm} \setlength{\parskip}{0cm}
\bproof
Let $U$ be the basic open subset
\bglnoz
&& U(x+I;x_1+I_1,\dotsc,x_n+I_n) \\
&=& \menge{\chi \in \widehat{E}}{\chi((x+I) \times I\reg) = 1; \, \chi((x_1+I_1) \times I_1\reg) = \dotsc = \chi((x+I) \times I\reg) = 0},
\eglnoz
where $I_i \subseteq I \subsetneq R$. By Lemma~\ref{co=disjointunion}, it suffices to show that $U$ is $(K \rtimes K\reg,\cC \cO)$-paradoxical.

Set $J \defeq \bigcap_{i=1}^n I_i$. Obviously, $J \neq (0)$. Therefore $1 + J \nsubseteq R^*$. Otherwise, by \cite[Proposition~1.9]{AM}, $J \subseteq {\rm Jac}(R) = (0)$, which is a contradiction. So we can choose $a \in (1+J) \setminus R^*$, $a \neq 0$, as $R\reg \neq R^*$. Then $J \nsubseteq aR$. Otherwise, there would exist $r \in R$ with $a-1 = ar$, so that $a(1-r)=1$ contradicting $a \notin R^*$. Choose $b_1 \in J$, $\delta \in J \setminus aR$ and set $b_2 \defeq b_1 + \delta$.
\setlength{\parindent}{0.5cm} \setlength{\parskip}{0cm}

For $k = 1,2$, $(b_k,a).U = U(b_k + ax+aI;b_k + ax_1+aI_1,\dotsc,b_k + ax_n+aI_n)$. Since $b_k + ax + aI \subseteq ax + aI + J \subseteq a(x+I) + J \subseteq (1+J)(x+I) + J \subseteq x+I$, every $\chi \in (b_k,a).U$ satisfies $\chi((x+I) \times I\reg) = 1$. Moreover, we have $(x_i + I_i) \cap (b_k + ax + aI) = b_k + (x_i + I_i) \cap (ax + aI) = b_k + (ax_i + I_i) \cap (ax_i + aI)$ since $x_i - ax_i = (1-a)x_i \in J$.

We claim that $I_i \cap aI = aI_i$. \an{$\supseteq$} is clear. If $r \in I_i \cap aI$, then $r = as$ for some $s \in I$, so that $s = (1-a+a)s = (1-a)s + as \in J + I_i \subseteq I_i$. This proves \an{$\subseteq$}.

Hence $(x_i + I_i) \cap (b_k + ax + aI) = b_k + ax_i + aI_i$. Therefore, every $\chi \in (b_k,a).U$ satisfies $\chi((x_i+I_i) \times I_i\reg) = 0$. This shows that for $k = 1,2$, $(b_k,a).U \subseteq U$.

Also, $(b_1 + ax + aI) \cap (b_2 + ax + aI) = \emptyset$ as $b_2 - b_1 = \delta \notin aR$. Hence $(b_1,a).U \cap (b_2,a).U = \emptyset$.
\eproof
\setlength{\parindent}{0cm} \setlength{\parskip}{0cm}
This strengthens and explains \cite[Theorem~1.3]{Li3}.
\setlength{\parindent}{0cm} \setlength{\parskip}{0.25cm}

\subsection{Almost finite groupoids}
Apart from purely infinite groupoids, Matui also introduced almost finite ones in \cite{Ma1}. We would like to end with the following obervation concerning the relation between almost finiteness and amenability:
\bprop
Let $G$ be a group acting on a compact space $X$. Assume that there exists $x \in X$ with trivial stabilizer group, i.e., $G_x = \menge{g \in G}{g.x = x} = \gekl{e}$. If the transformation groupoid $G \ltimes X$ is almost finite, then $G$ is amenable.
\eprop
\setlength{\parindent}{0.5cm} \setlength{\parskip}{0cm}
\bproof
Let $E \subseteq G$ be a finite subset with $E = E^{-1}$. For every $\ve > 0$, we have to find a finite subset $F \subseteq G$ such that for every $s \in E$, $\abs{sF \triangle F} / \abs{F} < \ve$.

Set $C \defeq E \times X \subseteq G \ltimes X$. Let $x \in X$ be as above, i.e., $\menge{g \in G}{g.x = x} = \gekl{e}$. Since $G \ltimes X$ is almost finite, we can find an elementary subgroupoid $K$ of $G \ltimes X$ such that $\abs{CKx \setminus Kx} / \abs{Kx} < \frac{\ve}{2}$. Let $F = \menge{g \in G}{(g,x) \in K}$. $F$ is contained in the image of $K$ under the canonical projection $G \ltimes X \to G, \, (g,x) \ma g$. Therefore, $F$ is a finite subset of $G$. We have $Kx = \menge{(g,x)}{g \in F}$ and $CKx \setminus Kx = \menge{(g,x)}{g \in EF \setminus F} = \bigcup_{s \in E} \menge{(g,x)}{g \in sF \setminus F}$. Therefore, $\frac{\ve}{2} > \abs{CKx \setminus Kx} / \abs{Kx} \geq \abs{sF \setminus F} / \abs{F}$ for all $s \in E$. For $s \in E$, $F \setminus sF = F \setminus (F \cap sF) = s(s^{-1}F \setminus (s^{-1}F \cap F) = s(s^{-1}F \setminus F)$. Hence $\abs{F \setminus sF} = \abs{s^{-1}F \setminus F} < \frac{\ve}{2} \abs{F}$ by our computation above. (Note that by assumption, $E = E^{-1}$.)

Therefore, for every $s \in E$, $\abs{sF \triangle F} / \abs{F} = \abs{sF \setminus F} / \abs{F} + \abs{F \setminus sF} / \abs{F} < \ve$.
\eproof
\setlength{\parindent}{0cm} \setlength{\parskip}{0.5cm}

\end{document}